\documentclass[11pt]{amsart}

\usepackage{amsmath,amsthm, amscd, amssymb, amsfonts}
\usepackage[all]{xy}
\usepackage{verbatim}
\usepackage{color}
\usepackage{array}
\usepackage{multirow}
\usepackage{hhline}
\usepackage{enumitem}

\usepackage{mathtools} 
\usepackage[normalem]{ulem} 

\newtheorem{theorem}{Theorem}[section]
\newtheorem{lem}[theorem]{Lemma}

\newtheorem*{claim}{Claim}

\newtheorem*{taskone}{Task 1}
\newtheorem*{tasktwo}{Task 2}

\theoremstyle{definition}
\newtheorem{definition}[theorem]{Definition}
\newtheorem{exa}[theorem]{Example}

\theoremstyle{remark}
\newtheorem{rem}[theorem]{Remark}

\newcommand{\ydh}{{}^{H}_{H}\mathcal{YD}}
\newcommand{\ydk}{{}^{K}_{K}\mathcal{YD}}

\newcommand{\ig}{:=}

\def\pf{\begin{proof}}
\def\epf{\end{proof}}

\newcommand{\ydgdual}{{}_{\ku^G}^{\ku^G}\mathcal{YD}}
\newcommand{\nc}{\newcommand}

\def\bs{\boldsymbol}

\newcommand{\I}{{\mathbb I}}
\newcommand{\Jb}{{\mathbb J}}

\nc{\ub}{\mathfrak{u}}
\nc{\g}{\mathfrak{g}}

\newcommand{\G}{{\mathbb G}}
\def\Bt{\widetilde{\B}}
\def\cHt{\widetilde{\mathcal{H}}}

\newcommand\GL{\operatorname{GL}}
\newcommand\id{\operatorname{id}}

\newcommand{\ydg}{{}^{\ku G}_{\ku G}\mathcal{YD}}

\newcommand\mo{\text{-}\operatorname{mod}}
\newcommand\sgn{\operatorname{sgn}}

\newcommand\gr{\operatorname{gr}}

\newcommand\co{\operatorname{co}}
\newcommand\ad{\operatorname{ad}}

\newcommand\rg{\operatorname{rg}}

\newcommand\ord{\operatorname{ord}}
\newcommand\Alg{\operatorname{Alg}}
\newcommand\Cleft{\operatorname{Cleft}}

\newcommand\Rep{\operatorname{Rep}}

\def\k{\Bbbk}
\def\ku{\Bbbk}

\def\ot{\otimes}

\def\Dm{\mathbb{D}}
\def\s{\mathbb{S}}
\def\C{\mathbb{C}}

\def\N{\mathbb{N}}
\def\Z{\mathbb{Z}}

\def\mH{\mathcal{H}}
\def\mB{\mathcal{B}}
\def\mA{\mathcal{A}}
\def\mG{\mathcal{G}}
\def\mJ{\mathcal{J}}

\def\Ss{\mathcal{S}}
\def\mL{\mathcal{L}}
\def\mP{\mathcal{P}}

\def\mO{\mathcal{O}}
\def\mE{\mathcal{E}}

\def\mH{\mathcal{H}}

\def\eps{\varepsilon}

\def\lg{\langle}
\def\rg{\rangle}

\newcommand{\Gc}{{\mathcal G}}

\newcommand{\qb}{{\bf q}}

\def\uDelta{\underline{\Delta}}

\renewcommand{\a}[1]{^{#1}}

\def\bx{\mathbf{x}}

\def\n{\mathfrak{n}}
\def\z{\mathfrak{z}}
\def\B{\mathfrak{B}}
\def\F{\mathfrak{F}}
\def\L{\mathfrak{L}}

\begin{document}



\title[How to lift]{Pointed Hopf algebras: a guided tour to the liftings}

\author[Angiono, Garc\'ia Iglesias]{Iv\'an Angiono, Agust\'in Garc\'ia Iglesias}

\address{FaMAF-CIEM (CONICET), Universidad Nacional de C\'ordoba,
Medina A\-llen\-de s/n, Ciudad Universitaria, 5000 C\' ordoba, Rep\'ublica Argentina.}

\email{(angiono|aigarcia)@famaf.unc.edu.ar}

\thanks{\noindent 2010 \emph{Mathematics Subject Classification.}
16T05. \newline The work was partially supported by CONICET,
FONCyT-ANPCyT, Secyt (UNC), the MathAmSud project GR2HOPF}

\begin{abstract}
This article serves a two-fold purpose. On the one hand, it is a survey about the classification of finite-dimensional pointed Hopf algebras with abelian coradical, whose final step is the computation of the liftings or deformations of graded Hopf algebras. On the other, we present a  step-by-step guide to carry out the strategy developed to construct the liftings. As an example, we conclude the work with the classification of pointed Hopf algebras of type $B_2$.
\end{abstract}

\maketitle

\section{Introduction}\label{sec:intro}

A pointed Hopf algebra $A$ is characterized by the fact that its coradical $A_0$ coincides with the subalgebra $\k G$ generated by its group-like elements $G=G(A)$. When considering the classification problem, this group $G$ is thus a first invariant. Associated to it there is also a braided structure at the heart of $A$: the so-called infinitesimal braiding; this is an object $V$ in the category of Yetter-Drinfeld modules $\ydg$. Starting with $G$ and $V$, the classification of finite-dimensional Hopf algebras with abelian coradical has been achieved throughout the collaborative work of many authors, specially those of Andruskiewitsch and Schneider \cite{AS-lift-meth, AS-annals}, Heckenberger \cite{H} and Angiono \cite{A-jems, Ang-crelle}. The final step was recently completed by the authors in \cite{AnG}, based on a strategy to construct Hopf algebras -the liftings- developed in \cite{AAGMV} and \cite{AAG}. In this article, we shall review the development of this classification and give detailed instructions about how to carry on this program on each example. This step-by-step guide is the self-contained Section \ref{sec:howtolift} and can be extracted by a potential user.

\subsection{The Diagonal Setting}\label{sec:settings}

One of the main purposes of this article is to present the recipe  to construct liftings in full detail; this is presented in \S \ref{sec:howtolift} and an example is developed in \S \ref{sec:example}. Although the strategy developed in \cite{AAGMV} applies in a more general level, we restrict ourselves for pedagogical reasons to the following setting: 
\begin{enumerate}
\item A cosemisimple Hopf algebra $H$.
\item A braided vector space of diagonal type such that $\dim \B(V)<\infty$ with a principal realization $V\in\ydh$.
\end{enumerate}
The strategy provides a recursive algorithm to construct liftings of $V$, that is Hopf algebras $A$ with $\gr A\simeq \B(V)\# H$.
It is important to notice that we do not ask $H$ to be finite-dimensional, nor we assume $H$ to be commutative. 

See \S \ref{sec:bvs} for unexplained notation, in particular Definition \ref{def:realization}.

\subsection{The lifting method}
The program for the classification of finite dimensional pointed Hopf algebras (over $\C$) with abelian group of group-like elements was originated by Andruskiewitsch and Schneider in \cite{AS-lift-meth}.  This class of Hopf algebras contains the small quantum groups $\ub^+_q(\g)$, $\g$ a semisimple Lie algebra and $q\neq 1$ a root of 1 in $\k$. 
Let $\Gamma$ be a finite abelian group and let $H$ be a pointed Hopf algebra with $H_0\simeq \k\Gamma$. A key observation that triggered the ulterior development was the fact that for each such $H$ there is a braided vector space $(V,c)$ of diagonal type and hence a graded algebra $\B(V)$ generated by the degree one component $V=R_{(1)}$ of the coinvariant subalgebra $R=\oplus_{k\geq 0} R_{(k)}=(\gr H)^{\co H_0}$ induced by the projection $\gr H\simeq R\# \k\Gamma\twoheadrightarrow \k\Gamma$. This can be generalized for any cosemisimple Hopf algebra $H_0$ such that it is a Hopf subalgebra of $H$.

\subsubsection{} The lifting method developed by Andruskiewitsch and Schneider to classify finite-dimensional pointed Hopf algebras $A$ with abelian coradical $A_0\simeq \k\Gamma$ consists of the several steps. Here we describe it for a general cosemisimple coradical $A_0\simeq K\subset A$ that is a Hopf subalgebra:
\begin{enumerate}
\item to classify all braided vector spaces $(V,c)$ with a principal realization in $\ydk$ such that $\dim\B(V)<\infty$.
\item to give a presentation of $\B(V)$.
\item to check if, given $H$ with $H_0\simeq K$, then $(\gr H)^{\co H_0}\simeq \B(V)$ (generation in degree one).
\item to compute all Hopf algebra deformations of $\B(V)\# K$ (to lift the relations of $\B(V)$).
\end{enumerate}

When $K=\k\Gamma$, any $V\in\ydk$ is necessarily equipped with a diagonal braiding, that is it has a linear basis $\{x_1,\dots,x_\theta\}$ and there is a matrix $\qb=(q_{ij})_{i,j\in\I_\theta}$ such that the braiding $c=c^{\qb}$ is determined by the equation:
\[
c(x_i\ot x_j)=q_{ij}\, x_j\ot x_i, \quad i,j\in\I_\theta.
\] 
In this setting, step (1) was completed by Heckenberger in \cite{H}, step (2) was achieved by Angiono in \cite{A-jems, Ang-crelle}, who used this result to prove (3) in \cite{Ang-crelle}. The main contribution of \cite{AAGMV, AAG,AnG}, which is the main focus of this survey, was to provide a strategy to complete (4) and to prove that it provided a full classification. This strategy is built on cocycle deformations of graded Hopf algebras, as suggested by a result of Masuoka \cite{M}, who showed that the class of liftings given in \cite{AS-annals} were cocycle deformations of the associated graded algebras. This phenomenon had been previously glimpsed in \cite{GM}. A different approach to solve step (4) in some cases in rank 2 was developed in \cite{Helbig}.

\subsubsection{} 
The strategy in \cite{AAGMV} consists in constructing a collection of liftings $\ub(\bs\lambda)$ indexed by a family of parameters $\bs\lambda\in\bs\Lambda$. These Hopf algebras are obtained as the final step in a sequence $(\mL_k=\mL_k(\bs\lambda))_k$ of Hopf algebra quotients 
\[
\mL_0=T(V)\# H\twoheadrightarrow \mL_1\twoheadrightarrow\dots  \twoheadrightarrow \mL_{\ell+1}=\ub(\bs\lambda).
\]

In turn, each $\mL_k$ is recovered as the Schauenburg left Hopf algebra associated to a cleft object $\mA_k(\bs\lambda)$ for a graded Hopf algebra $\mH_k=\B_k\# H$; here $(\B_k)_k$ is a suitably chosen family of pre-Nichols algebras so that
\[
\B_0=T(V)\twoheadrightarrow \B_1\twoheadrightarrow\dots  \twoheadrightarrow \B_{\ell+1}=\B(V).
\]
is a chain of braided Hopf algebra quotients in $\ydh$.

A remarkable fact is that each cleft object $\mA_k(\bs\lambda)$ splits as $\mA_k(\bs\lambda)=\mE_k(\bs\lambda)\# H$, for a certain coinvariant $H$-module algebra $\mE_k(\bs\lambda)$, $k\in\I_{\ell+1}$.

This was adapted for the diagonal setting in \cite{AAG} and showed to be an exhaustive method to complete the classification provided that $\mE(\bs\lambda)\coloneqq \mE_{\ell+1}(\bs\lambda)$ was nonzero; a fact that was ultimately proved in \cite{AnG}. Furthermore, in \cite[\S 5.2]{AnG} we provide a list of tools and techniques to deal with this problem in different settings.

\subsection{A non diagonal setting}\label{sec:more-general}
 The strategy we shall describe actually works in a more general framework, where we do not assume $V$ to be of diagonal type nor $\B(V)$ to be finite dimensional. More precisely, 
 it can be used to compute liftings of $V$ over $H$ when we have the following setup:
\begin{enumerate}
	\item A cosemisimple Hopf algebra $H$.
	\item A braided vector space such that the ideal $\mJ(V)$ defining the Nichols algebra $\B(V)$ is finitely generated.
\end{enumerate}
Indeed, this strategy is based on a recursive argument indexed by a minimal set $\Gc$ of generators of the ideal  $\mJ(V)$; hence we need not ask $\B(V)$ to be finite-dimensional, but rather that  $\mJ(V)$ is finitely generated. We invite the reader interested in this wider scope to check Section \ref{sec:gral}.  The cosemisimplicity hypothesis  on $H$ is, on the other hand, necessary: we will discuss this in \S \ref{sec:coss}. Finally, we discuss in \S \ref{sec:non-principal}  when we can remove the request for the realization $V\in\ydh$ to be principal.

\subsection{} The paper is organized as follows: in Section \ref{sec:pre} we write down all the preliminaries and notations needed throughout the paper, together with a brief survey on the development of the classification problem. We illustrate this with a toy example in \ref{sec:toy}.  In Section \ref{sec:strategy} we describe the strategy to compute the liftings. Section \ref{sec:howtolift} is a self-contained do-it-yourself operation manual to compute liftings of braided vector spaces of diagonal type; which we complete with an example in Section \ref{sec:example}, where we present the classification of liftings of diagonal braidings of type $B_2$. Finally, Section \ref{sec:gral} is devoted to the exploration of the outer limits of the strategy.

\subsection*{Acknowledgments.} This survey started with the suggestion of Nicol\'as Andruskiewitsch to the second author during the CLA in Quito, Ecuador. Both authors wish to thank him for his constant support and generosity along these years; most of the results in this note have been achieved in collaboration with him or motivated by him.

\section{Preliminaries}\label{sec:pre}

We work over an algebraically closed field $\k$ of characteristic
zero. For $\theta\in\N$, we set $\I_\theta\coloneqq\{1,\dots,\theta\}\subset \N$.
All algebras, tensor products, etc.~are considered over $\k$.
If $G$ is a group, we denote by $Z(G)$ the center of $G$. 

We write $\G_N$ for the group of $N$th roots of 1 in $\k$ and $\G_N'\subset \G_N$ for the subset of primitive roots.
We use $\GL(V)$ to denote the group of invertible linear maps on a vector space $V$. We denote the symmetric group on $n$ letters by $\s_n$ and use $\Dm_m$ to denote the dihedral group of order $2m$.

Let $A,B$ be algebras, we write $\Alg(A,B)$ for the set of algebra maps $A\to B$.
If $S\subset A$ is a set, we denote by $\lg S\rg\subseteq A$ the ideal generated by $S$.

If $H$ is a Hopf algebra with comultipication $\Delta$, then we
shall use Sweedler's notation $\Delta(h)=h_{(1)}\ot h_{(2)}$;
similarly for a (right) comodule $(M,\rho)$ over $H$: $\rho(m)=m_{(0)}\ot m_{(1)}$, $m\in M$.
We denote by  $G(H)=\{x\neq 0:\Delta(x)=x\ot x\}$ 
the group of grouplike elements of $H$. We write $(H_n)_{n\geq 0}$ for the
coradical filtration of $H$; namely $H_0=\sum_{C} C$ and
$H_{n+1}=\Delta^{-1}(H\ot H_n+H_0\ot H)$; where $C$ runs over all
simple subcoalgebras of $H$ -in particular $\k G(H)\subset H_0$. 

When $H_0$ is a Hopf subalgebra of $H$, the associated
graded coalgebra $\gr H=\oplus_{n\geq 0} H_{n+1}/H_n$ is actually 
Hopf algebra. A distinguished family of
Hopf algebras satisfying this condition is the class of \emph{pointed} Hopf algebras, which are
defined by requesting $H_0=\k G(H)$. Also, \emph{copointed} Hopf algebras (those with $H_0=\k^G$, $G$ a non-abelian group) satisfy this.

\subsection{Cocycle deformations and cleft objects}

\subsubsection{Cocycles}
A convolution-invertible linear map $\sigma: H\ot  H \to \ku$ on a Hopf algebra $H$ is said to be a 2-cocycle if the following holds:
\begin{align*}
\sigma(x_{(1)}, y_{(1)}) \sigma(x_{(2)} y_{(2)}, z) &=
\sigma(y_{(1)}, z_{(1)}) \sigma(x, y_{(2)}z_{(2)}), \ \text{all }x,y,z\in H.
\end{align*}
It can be normalized by setting: $\sigma(x, 1) = \sigma(1, x) = \varepsilon(x)$, $x\in H$.
We write $Z^2(H,\k)$ for the set of (normalized) 2-cocycles on $H$. Given $\sigma\in Z^2(H,\k)$, we may twist the multiplication $H\ot H\to H$ by setting:
\begin{align*}
x\cdot_{\sigma}y &= \sigma(x_{(1)}, y_{(1)})
 x_{(2)} y_{(2)}
\sigma^{-1}(x_{(3)}, y_{(3)}), \quad x, y\in H.
\end{align*}
This defines a new associative product on the vector space $H$ in such a way that $(H,\cdot_\sigma, 1, \Delta,\varepsilon)$ is again a Hopf
algebra, with a certain antipode $\Ss_{\sigma}$. 

\smallbreak

We denote this new Hopf algebra $H_\sigma$.

\begin{definition}
Let $A$, $H$ be Hopf algebras. Then $A$ is a cocycle deformation of $H$ if there is $\sigma\in Z^2(H,\k)$ such that $A\simeq H_\sigma$ as Hopf algebras.
\end{definition}

\subsubsection{Cleft objects}

Let  $C$ be a (right) $H$-comodule algebra: then $C$ is called a (right) cleft object for $H$ if $C^{\co H}=\k$ and there exists a convolution-invertible comodule isomorphism
$\gamma:H\to C$. We may further assume that $\gamma(1)=1$; such $\gamma$ is
called a {\it section}. We set 
\[
\Cleft (H):=\{\text{isomorphism classes of $H$-cleft objects}\}.
\]
Left, resp. bi-, cleft objects are defined analogously.

\subsubsection{G\"unther's approach}\label{sec:gunther}

In \cite{Gunther}, G\"unther develops a way to relate the cleft objects for a given Hopf algebra $H$ with the cleft objects $C'$ for a quotient Hopf algebra $H\twoheadrightarrow H'$; each $C'$ arises as a quotient $C\twoheadrightarrow C'$ of a certain $C\in\Cleft(H)$. Reciprocally, each $C\in\Cleft(H)$ can be recovered as a cotensor product $C\simeq C'\square H$ for some $C'\in\Cleft(H')$. As explained in \ref{sec:coss}, one needs to assume that $H$ is $H'$-coflat.

There are two alternatives to construct a quotient $C'$ (we \underline{underline} the grain of salt on each one):
\begin{enumerate}
\item Either we 
\begin{enumerate}
\item \underline{compute} $X=\,^{\co H'}H$,
\item and pick $\varphi\in\Alg_H^H(X,C)$;
\end{enumerate}
\item Or we
\begin{enumerate}
\item fix a right coideal subalgebra $Y\subset H$ such that $H'=H/\lg Y^+\rg$,
\item and pick $\varphi\in\Alg^H(Y,C)$ \underline{such that} $C\varphi(Y^+)C\neq C$. 
\end{enumerate}
\end{enumerate}
According to each alternative, we set
\begin{align}\label{eqn:cleft-quotient-gunther}
C'&=C/\lg\varphi(X^+)\rg & &\text{or else } & C'&=C/\lg\varphi(Y^+)\rg.
\end{align}

\subsubsection{Schauenburg's left Hopf algebra}
Let $H$ be a Hopf algebra and fix $C\in \Cleft(H)$. Then there is a Hopf algebra $L=L(C,H)$ in such a way that $C$ is a $(L,H)$-bicleft object. In particular, it follows that $L$ is a cocycle deformation of $H$. The converse is also true, if $A$ is a cocycle deformation of $H$, then there is $C\in \Cleft(H)$ such that $A\simeq L(C,H)$. See \cite{S} for details.

\subsection{Yetter-Drinfeld modules and Nichols algebras}\label{sec:nichols}
Let us fix $H$ a Hopf algebra with bijective antipode. The category
of (left) $H$ Yetter-Drinfeld modules $\ydh$ is that of simultaneously
(left) $H$-modules and (left) $H$-comodules $V$ for which the
following compatibility condition holds:
\begin{align*}
\lambda(h\cdot v)&=h_{(1)}v_{(-1)}\Ss(h_{(3)})\ot h_{(2)}\cdot
v_{(0)}, & & h\in H, v\in V.
\end{align*}
This is a braided tensor category, with braiding $c_{V,W}\colon V\ot W\to W\ot V$,
\begin{align*}
c_{V,W}(v\ot w)&=v_{(-1)}\cdot w\ot v_{(0)}, & & v\in V, w\in W.
\end{align*}

\begin{exa}
Let $H=\k G$ be the group algebra of a finite group $G$. Then $\ydh$ is the category of $G$-graded $G$-modules $V=\oplus_{g\in G} V_g$ such that 
\begin{align*}
g\cdot V_h &\subseteq V_{ghg^{-1}} & \text{for all } &g,h\in G. 
\end{align*}

In particular, if $G$ is abelian, then $\ydh$ is the category of $G$-graded $G$-modules with stable homogeneous components.
\end{exa}

The \emph{Nichols algebra} $\B(V)$ of $V$ is defined as the maximal graded braided Hopf
algebra quotient $T(V)\twoheadrightarrow \B$ such that $P(\B)=V$. We denote by $\mJ(V)\subset T(V)$ the ideal such that $\B(V)=T(V)/\mJ(V)$. When this is finitely generated, we write $\mG(V)\subset \mJ(V)$ for a minimal set of generators.

\begin{definition}\label{def:lifting}
A lifting of $V\in \ydh$ is a Hopf algebra $A$ such that $\gr
A\simeq \B(V)\# H$.
\end{definition}

\subsection{Braided vector spaces}\label{sec:bvs}

Recall that a braided vector space is a pair $(V,c)$ where $V$ is a vector space and $c\in\GL(V\ot V)$ is a solution to the braid equation:
\begin{equation}\label{eqn:braid}
(c\ot \id)(\id\ot\, c)(c\ot \id)=(\id\ot\, c)(c\ot \id)(\id\ot\, c).
\end{equation}
 In particular, any object in $V\in \ydh$ is a braided vector space with $c=c_{V,V}$.

\begin{definition}\label{def:realization}
Let $H$ be a Hopf algebra, $(V,c)$ a braided vector space.
\begin{enumerate}
\item A realization of $(V,c)$ over $H$ is a 
structure of Yetter-Drinfeld $H$-module on $V$ in such a way that
$c$ coincides with the categorical braiding $c_{V,V}$ in $\ydh$.
\item A realization $V\in\ydh$ is called principal when there is a basis $\{x_i:i\in \I\}$ of $V$ and elements $g_i\in H$, $i\in\I$, such that the coaction is given by $x_i\mapsto g_i\ot x_i$, $i\in\I$; in particular $g_i\in G(H)$.
\end{enumerate}
\end{definition}

\subsubsection{} The Nichols algebra of $V\in\ydh$ as in \ref{sec:nichols} only depends on the braiding
of $V$; so it can be defined for any braided vector space $(V,c)$.

We may extend Definition \ref{def:lifting} to this setting by saying that
$A$ is a lifting of $(V,c)$ -over $H$, if there is a realization of
$V\in \ydh$ and $A$ is a lifting of $V$ as in Definition \ref{def:lifting}.

\subsection{Diagonal type}\label{sec:diag-type}

A  braided vector space $(V,c)$ is called of diagonal type if there is a linear basis $\{x_1,\dots,x_\theta\}$ of $V$ and a collection of scalars $(q_{ij})_{i,j\in \I_\theta}$ such that $c(x_i\ot x_j)=q_{ij}\,x_j\ot x_i$, $i,j\in\I_\theta$. We refer to $\qb=(q_{ij})$ as the braiding matrix.
In this setting, a principal realization $V\in \ydh$ amounts to the existence of a collection $(g_i)_{i\in \I_\theta}\in Z(G(H))$ and a family $(\chi_i)_{i\in \I_\theta}\in \Alg(H,\k)$ satisfying
\begin{align}\label{eqn:yd-pair}
\chi_i(h)g_i&=\chi_i(h_{(2)})h_{(1)}g_i\Ss(h_{(3)}), & h\in H, i\in\I_\theta.
\end{align}
A pair $(g,\chi)$ as in \eqref{eqn:yd-pair} is called a YD-pair.

\subsection{The lifting method}\label{sec:lifting-method}

Let us fix $G$ a finite abelian group. In \cite{AS-lift-meth}, Andruskiewitsch and Schneider proposed the following approach to find all finite-dimensional pointed Hopf algebras $A$ with $G(A)\simeq G$.

In particular, $A_0=\k G\subset A$ is a Hopf subalgebra. In this case, the coradical filtration $A_0\subset A_1\subset\dots $ is a Hopf algebra filtration and thus the graded coalgebra $\gr A=\oplus_{n\geq 0} A_n/A_{n-1}$ is a Hopf algebra. 

Now, $\gr A$ splits into a semidirect product (the bosonization) $\gr A=R\# \k G$, where $R$ is the coinvariant subalgebra with respect to the projection $\gr A\twoheadrightarrow A_0$. Moreover, $R=\oplus_{n\geq 0}R_{(n)}$ is a braided graded Hopf algebra in the category $\ydg$, with $R_{(0)}=\k$. 

Set $V\coloneqq R_{(1)}\simeq A_1/A_0$; this is called the infinitesimal braiding of $A$. Then $V\in\ydg$ (hence it is a braided vector space of diagonal type) and the subalgebra generated by $V$ is the Nichols algebra $\B(V)$. 

These observations lead the authors in \cite{AS-lift-meth} to propose the following Lifting Method to achieve the classification.

\subsubsection{Step 1: Find all $V\in\ydg$ such that $\dim\B(V)<\infty$}\label{sec:lm-step1}
The underlying braided vector space of such a $V$ is necessarily of diagonal type. The classification of all (connected) $(V,c^{\qb})$ of diagonal type -- equivalently of all braiding matrices $\qb=(q_{ij})_{i,j}$ -- with $\dim\B(V)<\infty$ was completed in \cite{H}, in terms of so-called generalized Dynkin diagrams. These diagrams are connected decorated graphs with vertices $\{1,\dots,\theta\}$, for $\theta=\dim V$, and there is an edge connecting $i$ and $j$ if and only if $\widetilde{q_{ij}}\ig q_{ij}q_{ji}\neq 1$. Each vertex $i$ is decorated with the label $q_{ii}$ while each edge is decorated with the value $\widetilde{q_{ij}}$.

Each diagram is associated to a Weyl grupoid $\mathcal{W}$ and to a root system $\Delta$. These data split the diagrams into different classes, such as Cartan type, Super type, Standard type, Modular type and Unidentified type.

\begin{exa} 
In \S\ref{sec:example}, we shall present the classification of all the liftings of braided vector spaces of Cartan type $B_2$; in this case $\qb=\begin{pmatrix}
q&q_{12}\\q_{21}&q^2
\end{pmatrix}$, where $q\in\G_N'$, $N\geq 3$, and $q_{12}q_{21}=q^{-2}$. The generalized Dynkin diagram is thus
\begin{align*}
\xymatrix{  \overset{\,\,q}{\underset{\ }{\circ}} \ar  @{-}[r]^{q^{-2}} &
\overset{q^2}{\underset{\ }{\circ}}}
\end{align*}
\end{exa}

With this information, when the group $G$ is given, it just remains to check which matrices $\qb=(q_{ij})_{i,j}$ can be realized as a Yetter-Drinfeld module over $G$: this is a simple verification which amounts to solving linear equations involving the characters of the group. That is, whether there are elements $g_1,\dots,g_\theta\in G$ and characters $\chi_1,\dots,\chi_\theta\in\widehat{G}$ satisfying $\chi_j(g_i)=q_{ij}$, $i,j\in\I_\theta$.

\begin{exa}
Let $(V,c)$ be a braided vector space of diagonal type Cartan $B_2$, with $q\in\G_6'$ (so $q^2\in\G_3$).  Let 
\[
G=\Z/6\Z\times \Z/3\Z=\lg t_1,t_2:t_1t_2=t_2t_1, t_1^6=1, t_2^3=1\rg
\]
 and let $\tau_i\in\widehat{G}$ be such that $\tau_i(g_j)=q^{\delta_{i,j}i\cdot j}$. It follows that $V\in\ydg$ with $g_i=t_i$, $i=1,2$, and $\chi_1=\tau_1\tau_2$, $\chi_2=\tau_1^a\tau_2^b$, with $a+2b=4$ (e.g.~$\chi_2=\tau_1^2\tau_2$).
\end{exa}

\subsubsection{Step 2: For each such $V$, describe $\B(V)$ by generators and relations}\label{sec:lm-step1prime}

Starting with the classification of arithmetic root systems in \cite{H}, in particular of all braided vector spaces of diagonal type $(V,c^{\qb})$ with $\dim\B(V)<\infty$, described in \ref{sec:lm-step1}, Angiono computed the presentation of each Nichols algebra $\B(V)$, generated by $V$. As a byproduct, it follows that $\B(V)$ has a PBW basis with generators indexed by the positive roots $\Delta^+$ associated to $\qb$. The generators $x_i$ correspond to the simple roots $\alpha_i$, $i\in\I_\theta$ and generators $x_\alpha$, $\alpha\in\Delta^+$ are defined recursively starting with $x_{\alpha_i}=x_i$.

The set of defining relations $\Gc(V)$ can be generically split into two distinguished subsets, namely
\begin{itemize}
\item (generalized) quantum Serre relations, like $(\ad_{c}x_i)^{m_{ij}+1}(x_j)$, and
\item powers of root vectors of Cartan type\footnote{In a few, specific, examples, the square of a non-Cartan root $\alpha_i+\alpha_{i+1}$ is needed. As well, a simple root may  not be of Cartan type, but the corresponding power is also needed.} $x_\alpha^{N_\alpha}$.
\end{itemize}
The first item guarantees the skew-commutation of the letters $x_\alpha$, $\alpha\in\Delta$, while the second ensures the finite height of those of Cartan type (or simple). This is enough to cap the height of non-Cartan root vectors.

\begin{exa}
Let $(V,c)$ of Cartan type $B_2$; the root system is
\[
\alpha_1,\alpha_1+\alpha_2,2\alpha_1+\alpha_2,\alpha_2
\] 
and a PBW basis for $\B(V)$ is given by the set
\[
\{x_2^ax_{12}^bx_{112}^cx_1^d:0\leq a,c < M,0\leq b,d <N\},
\]
where $M=N$ if $N$ is odd and $M=N/2$ otherwise. 

The set of relations $\Gc(V)$ is given by:
\begin{itemize}
\item $(\ad_{c}x_1)^3(x_2)$, $(\ad_{c}x_2)^2(x_1)$,
\item $x_2^M$, $x_{12}^N$, $x_{112}^M$, $x_1^N$.
\end{itemize}
\end{exa}

\subsubsection{Step 3: Decide if any $A$ with $G(A)\simeq G$ satisfies $R=\B(V)$} 

This step proposes to check if, for any $H$ with $H_0\simeq \k G$, the coinvariant subalgebra $R\subset \gr H$ coincides with the Nichols algebra $\B(V)$, that is if $R$ is itself generated in degree one.

This was first conjectured by Andruskiewitsch and Schneider in \cite{AS-annals}, where it was proved for a large subclass of diagrams. Later on, Angiono showed that the conjecture was valid for any diagram in the list of \cite{H}.

\subsubsection{Step 4: Compute all deformations of $\B(V)\#\k G$} 
That is, this step involves computing all Hopf algebras $A$ such that $\gr A\simeq \B(V)\# \k G$.

By the previous step, if $H$ is such that $H_0=\k G$, then $\gr H=\B(V)\# \k G$ for some $V$. That is, any $H$ is a deformation of $\B(V)\# \k G$ in the sense that the defining relations of $H$ descend to the defining relations of $\B(V)$ when we consider the filtration induced by the coradical. Hence the computation of all Hopf algebras $H$ with $H_0\simeq \k G$ is achieved by lifting the relations of the Nichols algebra. 

See \S \ref{sec:toy} for a simple example depicting the concept of lifting a given relation. 

\subsubsection{Step 5: Check if any such $A$ is a cocycle deformation of $\B(V)\#\k G$}
Actually, this step is not part of the original schema settled in \cite{AS-lift-meth} by Andruskiewitsch and Schneider. However, Masuoka in \cite{M} proved that every lifting computed with the method in \cite{AS-annals} was indeed a cocycle deformation of $\B(V)\# \k G$. This was also the case for examples arising from non-abelian groups. Hence  this new step became a natural question. 

\subsection{Merging (and solving) Steps 4 and 5}
As we explained in \S\ref{sec:lifting-method}, Steps 1, 2 and 3 of the Lifting Method were completely solved by Heckenberger and Angiono.

The heart of the idea originated in \cite{AAGMV} to compute liftings is to proceed with Steps 4 and 5 at the same time, constructing cocycle deformations which are liftings and showing that these are all the liftings. We defined a family of Hopf algebras $\ub(\bs\lambda)$ and we showed in \cite{AnG} that this list was exhaustive. 

We devote Section \ref{sec:strategy} to give a step-by-step guide on how to produce each lifting $\ub(\bs\lambda)$ explicitly.

\subsection{Cosemisimplicity}\label{sec:coss} Our strategy is built in the computation of cleft objects for a Hopf algebra $A'$ out of a given set of cleft objects for another Hopf algebra $A$ with $A\twoheadrightarrow A'$. This is based in the work of G\"unther, which requires that $A$ is $A'$-coflat. Now, we actually have $A=R\# H$, $A'=R'\# H$ and the surjection $A\twoheadrightarrow A'$ is induced by a (braided) Hopf algebra surjection $R\twoheadrightarrow R'$ in $\ydh$. Moreover, it follows that $R$ is left and right cofree over $R'$. Since we assume that $H$ is cosemisimple, it follows that the coextension
\[
\id\ot\,\eps :R'\#H\to H
\]
is cosemisimple and thus $A$ is left and right cofree (hence coflat) over $A'$.

A {\bf question} arises: what is the situation if $H$ is not cosemisimple but generated, as an algebra, by a cosemisimple subcoalgebra $C$? This is the setting of the generalized lifting method introduced by Andruskiewitsch and Cuadra in \cite{AC}.

\subsection{A toy example}\label{sec:toy}

Let $V=\k\{x\}$ be a one-dimensional braided vector space, with braiding $c(x\ot x)=q\,x\ot x$, $q\in\G_N'$. 
The Nichols algebra $\B(V)$ is a truncated polynomial algebra: the defining ideal $\mJ(V)$ is generated by
\[\Gc=\{x^N\}.\]

Let $H$ be a cosemisimple Hopf algebra with a YD-pair $(g,\chi)$ as in \eqref{eqn:yd-pair} such that $\chi(g)=q$; that is to say, there is a realization $V\in\ydh$. 

\begin{exa}
We may take $H=\k G$, $G=\Z/mN\Z$, in which case the bosonization $\B(V)\# H$ is the generalized Taft algebra.
\end{exa}

The liftings of $V$ are given by a one-parameter family of Hopf algebras $\ub(\lambda)$, $\lambda\in\C$, defined as the quotient of $\k[a]\# H$ by the ideal generated by 
\[
\Gc(\lambda)=\{a^N-\lambda(1-g^N)\}.
\]
Here $\lambda$ is subject to the condition $\lambda=0$ if $\chi\neq \eps$.

Notice that, if $g^N=1$ in $H$, then $\ub(\lambda)= \B(V)\# H$ for all $\lambda$ and there is a unique (trivial) lifting. On the other hand, if $g^N\neq 1$, then there are two isomorphism classes, namely the class of $\lambda=0$ (corresponding to the trivial lifting) and the class of  $\lambda=1$.

In any case, the relation $x^N=0$ is lifted to a relation $a^N=\lambda(1-g^N)$.

Each Hopf algebra $\ub(\lambda)$ is a cocycle deformation of $\B(V)\#H$ by considering the bicleft object $\mA(\lambda)=\mE(\lambda)\# H$, where $\mE(\lambda)=\k[ y]/\lg\Gc'(\lambda)\rg$, for
\[
\Gc'(\lambda)=\{y^N-\lambda\}.
\]

\section{The strategy}\label{sec:strategy}

\subsection{Introduction}

Fix $H$, $V\in\ydh$ as in \ref{sec:settings}. We briefly review the strategy developed in \cite{AAGMV} to compute the liftings of $V$. Later on, we shall focus on each step, giving detailed instructions to perform them. Recall that a lifting of $V$ is a Hopf algebra $L$ such that $\gr L\simeq \B(V)\# H$; we set $\mH=\B(V)\#H$.

The main  objective, which is
\begin{center}
to find all Hopf algebras $L$ with $\gr L\simeq \mH$
\end{center}
is translated into a, possibly, less comprehensive one:
\begin{center}
to find all cocycle deformations $\mH_\sigma$ of $\mH$ such that $\gr \mH_\sigma\simeq \mH$.
\end{center}

\begin{rem}
When $V$ is of diagonal type and $\dim \B(V)<\infty$, then the two objectives listed above (to find all liftings and to find all cocycle deformations that are liftings) are actually equivalent by \cite{AnG}. The same is true for every $V$ of non-diagonal type with $\dim \B(V)<\infty$ that has been studied in the literature.
\end{rem}

Now, finding all  $\mH_\sigma$, that is finding all $\sigma\in Z^2(\mH,\k)$, is equivalent to finding all cleft objects $A\in\Cleft(\mH)$, as we can recover $\mH_\sigma$ as the Schauenburg left Hopf algebra $L(\mA,\mH)$ for a given $\mA\in \Cleft(\mH)$. Hence our objective is furthermore translated into
\begin{center}
to find all cleft objects $\mA\in\Cleft(\mH)$ such that $\gr L(\mA,\mH)\simeq \mH$.
\end{center}
That is precisely where the strategy leads to.

\subsection{General procedure}
Let us write $\Cleft'(\mH)\subset \Cleft(\mH)$ for the subset of those (isoclasses of) cleft objects $\mA$ such that $\gr L(\mA,\mH)\simeq \mH$; the subset that interests us.

A small hint towards the characterization of $\Cleft'(\mH)$ is given in \cite[Proposition 5.8]{AAGMV}, which states that two necessary conditions for a cleft object $\mA$ to be in the subset $\Cleft'(\mH)$ are
\begin{itemize}
\item there is an algebra surjection $\tau:T(V)\# H\twoheadrightarrow \mA$,
\item there is a section $\gamma:\mH\to \mA$ such that $\gamma_{|H}\in\Alg(H,\mA)$.
\end{itemize}

The procedure to find $\Cleft'(\mH)$ is recursive, using the ideas in \ref{sec:gunther}. More precisely, we consider a chain of Hopf algebra quotients:
\begin{align}\label{eqn:chain}
T(V)\#H=\mH_0\twoheadrightarrow \mH_1\twoheadrightarrow \dots \twoheadrightarrow \mH_N=\mH
\end{align}
and we deduce $\Cleft'(\mH_{k+1})$ from $\Cleft(\mH_k)$; starting with $\Cleft'(\mH_0)=\{\mH_0\}$. In particular, we get a list of Hopf algebras $\mL_k=L(\mH_k,\mA_k)$ for each $k=0,\dots, N$ and each $\mA_k\in\Cleft'(\mH_k)$. Two distinguished features in this sequence are
\begin{itemize}
\item $\mL_0\simeq T(V)\# H$,
\item If $\mA_k\in \Cleft'(\mH_k)$ projects onto $\mA_{k+1}\in \Cleft'(\mH_{k+1})$, then there is a Hopf algebra projection $\mL(\mA_k,\mH_k)\twoheadrightarrow\mL(\mA_{k+1},\mH_{k+1})$.
\end{itemize}
In particular, we recover all the liftings
$(L(\mA,\mH))_{\mA\in\Cleft'(\mH)}$ as quotients of $T(V)\# H$. This
amounts to saying that the output of the strategy is a list of
liftings of $V$ presented as algebras generated by $V$ and $H$, with some
relations. That is, the liftings are constructed in an explicit way.

\begin{rem}
A note is worthy to be mentioned: if $L$ is a lifting of $V$, e.g.~the Hopf algebras $L=L(\mA,\mH)$ for $A\in \Cleft'(\mH)$, then the coradical filtration $(L_n)_{n\geq 0}$ of $L$ satisfies $L_0\simeq H$ and $L_1\simeq V\# H$. When dealing with $\mH_k$, $k<N$, then this is not the case for $L=L(\mA_k,\mH_k)$, as this Hopf algebra  will provide a strict inclusion $V\# H\subset L_1$. The graded object we consider is associated to the filtration $\F=(F_n)$ induced by the graduation of $\mH_k$: that is, if $\pi_k:T(V)\#H\to\mH_k$ is the natural algebra projection, then $F_n=\pi_k(T^n(V)\#H)$. That is
\[
\Cleft'(\mH_{k})=\{\mA_k\in\Cleft(\mH_k)|\gr_{\F}L(\mA_k,\mH_k)\simeq \mH_k \}, \ k\leq N.
\]
In particular, when $k=N$, then this filtration coincides with the coradical filtration of $L(\mA,\mH)$ and this does not affect our previous definition of $\Cleft'(\mH)$.
\end{rem}

\subsection{The stratification}
To produce a suitable chain of quotients as in \eqref{eqn:chain} we stratify the minimal set of generators $\Gc=\Gc(V)$ of the ideal $\mJ(V)$ such that $\B(V)=T(V)/\mJ(V)$, as computed in \cite[Theorem 3.1]{Ang-crelle}. 

This stratification 
\[
\Gc=\Gc_0\sqcup\Gc_1\sqcup\dots \sqcup\Gc_\ell
\]
is chosen so that (the image of) the stratum $\Gc_{k}$ is contained in the set of primitive elements of the braided Hopf algebra 
$\B_k\coloneqq T(V)  / \lg \cup_{j=0}^{k-1} \Gc_j  \rg$. Here $\B_0=T(V)$. We define $\mH_k=\B_k\# H$ and this determines a sequence of $\ell+1$ recursive steps that we shall describe in \ref{sec:recstep} next.

\subsection{Forget $H$}\label{sec:forgetH}

A key observation regarding the algebras $\mA_k\in\Cleft'(\mH_k) $ is the following:
\begin{itemize}
\item Each $\mA_k$ splits as a smash product $\mA_k=\mE_k\# H$, where $\mE_k$ is an $H$-module algebra. 
\end{itemize}
Moreover, the section $\gamma_k:\mH_k\to \mA_k$ restricts to a braided $\mB_k$-comodule isomorphism $\gamma_k\ig\gamma_{k|\B_k}:\B_k\to \mE_k$, see \cite{AG}.

\subsection{The recursive step}\label{sec:recstep}

In this part we describe the mechanics to produce the set $\Cleft'(\mH_{k+1})$ and the liftings $\mL_{k+1}$, starting with the setting at a given level $k$.
That is, we assume that we have already computed $\Cleft'(\mH_k)$ and the liftings $\mL_k$. Notice that this is immediate
for $k=0$, as $\Cleft'(\mH_0)=\{T(V)\# H\}$ and $\mL_0=T(V)\# H$.

\subsubsection{Setting} We fix $k$ and some related notation: let $n_k$ denote the cardinality of the set of relations $\Gc_k$ and let $\{r_1,\dots,r_{n}\}$, $n=n_k$, be an enumeration of its elements. 

Now, we have a collection of algebras $\mA_k\in\Cleft \mH_k$, with a section $\gamma_k:\mH_k\rightarrow\mA_k$, and Hopf algebras $\mL_k=\mL_k(\mA_k,\mH_k)$ coacting on $\mA_k$ on the left, via $\delta_k:\mA_k\to\mL_k\ot\mA_k$.

 Recall from \ref{sec:forgetH} that $\mA_k=\mE_k\# H$, for some $\mE_k\in H\mo$.

\subsubsection{Input}
We have:
\begin{itemize}
\item[I.] a collection of algebras 
\begin{align}\label{eqn:Ek}
\{\mE_k\in H\mo\},
\end{align}
each of them with a section $\gamma_k:\B_k\rightarrow\mE_k$.
\item[II.] For each $r_i\in\Gc_k$, there is a YD-pair $(g_i,\chi_i)\in G(H)\times\Alg(H,\k)$ such that the $H$-action and coaction are determined by
\begin{align}\label{eqn:act-coact}
&& h\cdot r_i&=\chi_i(h)r_i, & r_i&\mapsto g_i\ot r_i, \qquad h\in H, \ i\in\I_n.
\end{align}
By definition, we have $\Delta(r_i)=r_i\ot 1+1\ot r_i$ in $\B_k$.
\item[III.] For each $r_i\in\Gc_k$, let $\tilde{r_i}\in\mL_k$ be such that
\begin{align}\label{eqn:deltak}
\delta_k(\gamma_k(r))-g_i\ot \gamma_k(r_i)=\tilde{r_i}\ot 1;
\end{align}
such $\tilde{r_i}$ exists by \cite[Proposition 5.10]{AAGMV}
\item[IV.] Finally, we consider all sequences 
\begin{align}\label{eqn:Lambdak}
\bs\lambda\in \bs\Lambda_k\ig \{(\lambda_i)_{i\in\I_n}\in\k^n:\lambda_i=0\text{ if }\chi_i\neq \eps\}.
\end{align}
\end{itemize}

\subsubsection{Output}\label{sec:output}
With these data in hand, we set 
\begin{align}
\label{eqn:Ek-quotient}
\mE_{k+1}(\bs\lambda)&=\mE_k/\lg\gamma_k(r_i)-\lambda_i: i\in\I_n\rg\\
\label{eqn:Lk-quotient}
\mL_{k+1}(\bs\lambda)&=\mL_k/\lg\tilde{r_i}-\lambda_i(1-g_i): i\in\I_n\rg.
\end{align}
We refer the interested reader to \ref{sec:quotients} for details on this construction.
Now, the collection
\[
\{\mE_{k+1}(\bs\lambda): \mE_k\text{ as in \eqref{eqn:Ek}}, \bs\lambda\in\bs\Lambda_k \}
\]
is a family of algebras in $H\mo$, each of them equipped with a section $\gamma_{k+1}=\gamma_k(\bs\lambda):\B_{k+1}\rightarrow\mE_{k+1}(\bs\lambda)$. 

Moreover, if $\mA_{k+1}(\bs\lambda)=\mE_{k+1}(\bs\lambda)\#H$, then 
 \[
 \Cleft'(\mH_{k+1})=\{\mA_{k+1}(\bs\lambda): \mE_k\in \Cleft'(\mH_k), \bs\lambda\in\bs\Lambda_k\}_{ / \sim}
  \]
  and $L(\mA_{k+1},\mH_{k+1})\simeq \mL_{k+1}(\bs\lambda)$.
%

\subsection{Final Output: the liftings}

Now we can compute the liftings of $V$ over $H$ as the Hopf algebras $\mL_{\ell+1}(\bs\lambda)$ arising in the last step of the strategy.

\begin{definition}\label{def:Lambda}
The set of deformation parameters is 
\begin{equation}\label{eqn:Lambda}
\bs\Lambda=\{(\lambda_r)_{r\in\Gc}:\lambda_r=0\text{ if }\chi_r\neq \eps\}.
\end{equation}
\end{definition}
It follows that, starting with $H$ and $V$, we obtain, for each $\bs\lambda\in\bs\Lambda$:
\begin{itemize} 
\item A family of $H$-module algebras $\mE(\bs\lambda)=T(V)/\lg \Gc'(\bs\lambda)\rg$, where 
\[
\Gc'(\bs\lambda)\ig\{\gamma_k(r)-\lambda_r: r\in\Gc_k, 0\leq k\leq\ell\}.
\]
\item A family of Hopf algebras $\ub(\bs\lambda)=T(V)\# H/\lg \Gc(\bs\lambda)\rg$, where 
\[
\Gc(\bs\lambda)\ig\{\tilde{r}-\lambda_r(1-g_r): r\in\Gc\}.
\]
\end{itemize} 

By \cite[Proposition 3.8]{AnG}, $\mE(\bs\lambda)\neq 0$ for every $\bs\lambda\in\bs\Lambda$; hence the family $(\ub(\bs\lambda))_{\bs\lambda\in\bs\Lambda}$ is the answer to the lifting problem:

\begin{theorem}Let $V$ be a braided vector space of diagonal type, with a realization $V\in\ydh$. Let $\bs\Lambda$ be as in \eqref{eqn:Lambda}.
For each $\bs\lambda\in\bs\Lambda$, the Hopf algebras $\ub(\bs\lambda)$ satisfy
\begin{enumerate}[leftmargin=*]
\item\cite{AAGMV} They are presented by generators and relations, as quotients of $T(V)\#H$, and obtained as the final step in a recursively defined chain of quotients.
\item\cite{AAGMV} $\gr \ub(\bs\lambda)\simeq \B(V)\# H$; i.e.~$\ub(\bs\lambda)$ is a lifting of $V$ over $H$.
\item\cite{AAG} $\ub(\bs\lambda)$ is a cocycle deformation of $\B(V)\# H$.
\end{enumerate}
Moreover, 
\begin{enumerate}[leftmargin=*]
\item[(4)]\cite{AnG} If $L$ is a lifting of $V$, then there is $\bs\lambda\in\bs\Lambda$ such that $L\simeq \ub(\bs\lambda)$.\qed
\end{enumerate}
\end{theorem}

\begin{rem}\label{rem:isos}
Two Hopf algebras $\ub(\bs\lambda), \ub(\bs\lambda')$, with $\bs\lambda\neq\bs\lambda'$, may be isomorphic. The isomorphism classes of these liftings are completely determined in \cite[Section 5]{AAG} in terms of an equivalence relation $\bs\lambda\sim \bs\lambda'$ in $\bs\Lambda$. Namely, each lifting is univocally determined by a single class $[\bs\lambda]$ in $\bs\Lambda/\sim$.
\end{rem}

\begin{rem}\label{rem:connected}
The classification of all $(V,c)$ of diagonal type, that is of all braiding matrices $\qb$, with $\dim\B(V)<\infty$ in \cite{H} is in terms of connected generalized Dynkin diagrams. If the diagram associated to a given $(V,c)$ is not connected, then the liftings of $V$ are determined by the liftings of its connected components, up to linking relations. That is, the lifting problem reduces to the connected case, see  \cite[Section 4]{AAG}.
\end{rem}

The same ideas as those discussed in Remark \ref{rem:connected} can be used to deal with subdiagrams. 
Next, we include this result for completeness.

Let $\I=\I_\theta$; we fix a basis $\{x_i\}_{i\in\I}$ of $V$ as in \S \ref{sec:diag-type} and let $\Jb\subset \I$. We set $W\subset V$ the braided subspace spanned by $\{x_j\}_{j\in\Jb}$. We shall relate the liftings of $W$ with certain subalgebras of liftings of $V$; notice that a realization $V\in\ydh$ immediately restricts to a realization $W\in\ydh$.

The Nichols algebra $\B(W)$ is a subalgebra of $\B(V)$, \cite{AS-pointed}. Moreover, $\Gc \cap T(W)\subset \Gc$ is the minimal set of defining relations 
$\Gc(W)$ from \cite{Ang-crelle}. As well, a stratification $\Gc=\Gc_0\sqcup\Gc_1\sqcup\dots \sqcup\Gc_\ell$ determines a stratification $\Gc(W)=\bigsqcup_k \Gc_k(W)$; with $\Gc_k(W)\coloneqq\Gc_k\cap T(W)$, $k=0,\dots,\ell$ (possibly empty for some $k$).

Now, given $\bs\lambda=(\lambda_r)_{r\in \Gc}\in\bs\Lambda$, we write $\bs\lambda_{|\Jb}\coloneqq (\lambda_r)_{r\in\Gc(W)}$. Notice that $\bs\lambda_{|\Jb}\in\bs\Lambda(W)$ and thus $\ub(\bs\lambda_{|\Jb})$ is a lifting of $W$. 
Here we denote by $\bs\Lambda(W)$ the set of deforming parameters corresponding to $H$ and $W$.

\begin{lem}\label{lem:subalgebra}
The Hopf algebra $\ub(\bs\lambda_{|\Jb})$ is (isomorphic to) the subalgebra of $\ub(\bs\lambda)$ generated by $H$ and $\{x_j\}_{j\in\Jb}$.
\end{lem}

\pf
Follows as the proof of \cite[Lemma 4.1]{AnG}, mutatis mutandis.
\epf

\subsection{On the quotients \eqref{eqn:Ek-quotient} and \eqref{eqn:Lk-quotient}}\label{sec:quotients}

Notice that, in other words, \eqref{eqn:act-coact} says that $\k\{r_i\}\subset \Gc_k$ is a sub-object in $\ydh$, for each $i\in\I_n$. 
%
%
%
%
In particular, 
\begin{align}\label{eqn:assume}
\text{we may assume $\Gc_k=\{r\}$ for a single $r$; }
\end{align}
we write $g_r\in G(H)$, $\chi_r\in\Alg(H,\k)$ for the corresponding structural data in \eqref{eqn:act-coact} and we set $q\ig \chi_r(g_r)$. If $Y=Y_k$ is the subalgebra of $\mH_k$ generated by $s\ig -\Ss(r)=rg_r^{-1}$, then it follows that
\begin{equation}\label{eqn:Yk}
Y
\simeq
\begin{cases}
\k[t], & \text{ if } q=1;\\
\k[t]/\lg t^N\rg, & \text{ if } q\neq 1 \text{ and }N=\ord q.
\end{cases}
\end{equation}
Also, $\mH_{k+1}=\mH_k/\lg r\rg=\mH_k/\lg s\rg=\mH_k/\lg Y^+\rg$, cf.~\ref{sec:gunther}.

Consider the (right) comodule map $\varphi=\varphi_{\lambda}:\k\{s\}\to \mA_k$, 
\begin{equation}\label{eqn:varphi}
\varphi(s)=\gamma_k(s)-\lambda g^{-1}=(\gamma_k(r)-\lambda) g^{-1}. 
\end{equation}
\begin{claim}
$\varphi$ extends to an algebra map $\varphi\colon Y\to \mA_k$.
\end{claim}
Hence, if  $\mA_k'=\mA_k/\lg \varphi(s)\rg\simeq \mA_k/\lg \gamma_k(r)-\lambda \rg$ is nonzero then it is a $\mH_{k+1}$-cleft object, see \eqref{eqn:cleft-quotient-gunther}. 

We check the claim: we may assume $Y\simeq \k[t]/\lg t^N\rg$; that is $r^N=0$. In particular, $\chi_r\neq \eps$. 
We need to show that $\varphi(s)^N=0$ or, equivalently, that $\gamma_k(r)^N=0$.
Now,
\begin{align*}
\rho(\gamma_k(r)^N)&=\rho(\gamma_k(r))^N=(\gamma_k(r)\ot 1+g_r\ot r)^N\\
&=\sum_{j=0}^N{ \binom{N}{j}_q}\gamma_k(r)^{N-j}g_r^j\ot r^j=\gamma_k(r)^N\ot 1+g_r^N\ot r^N\\
&=\gamma_k(r)^N\ot 1.
\end{align*}
Thus $\gamma_k(r)^N\in\mA^{\co \mH}\simeq \k$; set $c\ig \gamma_k(r)^N$. Let $\pi:\mA_k\to \mA_k'$ denote the algebra projection. We have that
$c=\pi(c)=\pi(\gamma_k(r)^N)=\pi(\gamma_k(r))^N=0$ and the claim follows.

\section{How to lift}\label{sec:howtolift}
{\it This section can be extracted from the article and used as a quick guide to compute liftings of diagonal type.}

\

Let $V,H$ as in \ref{sec:settings}. That is, $V=\k\{x_1,\dots,x_\theta\}$ is a braided vector space of diagonal type with $\dim\B(V)<\infty$ and $H$ is a cosemisimple Hopf algebra with a principal realization $V\in\ydh$. 

\subsection{Setting}

There are $(g_i)_{i\in\I_\theta}\in G(H)$ and $(\chi_i)_{i\in\I_\theta}\in\Alg(H,\k)$ such that the $H$-action and coaction on $V$ are given via $h\cdot x_i=\chi_i(h)x_i$ and $x_i\mapsto g_i\ot x_i$ and $\chi_j(g_i)=q_{ij}$, $i,j\in\I_\theta$. Here, $\qb=(q_{ij})_{i,j\in\I_\theta}$ stands for the braiding matrix of $V$, so that $c(x_i\ot x_j)=q_{ij}\,x_j\ot x_i$.

If $\Gc$ is a minimal set of homogeneous relations for $\B(V)$, then this setting determines elements $(g_r)_{r\in\Gc}\in G(H)$ and $(\chi_r)_{r\in\Gc}\in\Alg(H,\k)$.

Set $\mH=\B(V)\# H$ and fix a stratification $\Gc=\Gc_0\sqcup \Gc_1\sqcup \dots \sqcup \Gc_\ell$. This choice determines a recursive sequence of steps $0,\dots,\ell$. 

We write $\Gc^k=\cup_{j=0}^{k}\Gc_k$ and set $\B_{k+1}=T(V)/\lg \Gc^k\rg$, $\mH_{k+1}=\B_{k+1}\# H$; also $\B_0\ig T(V)$ so $\mH_0=T(V)\# H$.

In this setting, $\Delta(r)=r\ot1+1\ot r$ in $\B_k$, for $r\in \Gc_k$.

\subsection{A useful remark}\label{sec:remark-rho-delta}
All algebras involved in this procedure are generated by $V$ (and $H$). At any given step $k$, we denote the basis of $V$ by $\{x_i\}_{i\in\I_\theta},\{y_i\}_{i\in\I_\theta}$ or $\{a_i\}_{i\in\I_\theta}$ when we are dealing with $\B_k$, $\mA_k$ or $\mL_k$, respectively. We make no distinction on $H$ as it is naturally a subalgebra of all three algebras.

The coactions $\rho\colon\mA_k\to \mA_k\ot \mH_k$ and $\delta\colon\mA_k\to \mL_k\ot \mA_k$ are determined by $\Delta$ in the following way:
\begin{align}
\label{eqn:rho-i}
\rho(y_i)&=y_i\ot 1+g_i\ot x_i, \ i\in\I_\theta; & \rho(h)=h_{(1)}\ot h_{(2)}, \ h\in H;\\
\label{eqn:delta-i}\delta(y_i)&=a_i\ot 1+g_i\ot y_i, \ i\in\I_\theta; & \delta(h)=h_{(1)}\ot h_{(2)}, \ h\in H.
\end{align}

\subsection{The cleft objects}

We first look for the collection of algebras $\mA\in\Cleft'(\mH)$. Each one of them is determined by an $H$-module algebra $\mE=\mE(\bs\lambda)$, depending on a collection of scalars $\bs\lambda=(\lambda_r)_{r\in\Gc}$ in 
\[
\bs\Lambda=\{(\lambda_r)_{r\in\Gc} : \lambda_r=0 \quad  \text{ if } \quad \chi_r\neq\eps \}.
\]
in such a way that $\mA=\mE\# H$. Moreover, we have a precise presentation of each $\mE$, namely:
\[
\mE=T(V)/\lg \gamma_k(r)-\lambda_r: r\in \Gc_k \rg.
\]
So what is actually missing is to determine the elements
\[
\gamma_k(r)\in \mA_k, \ r\in \Gc_k.
\]
This will be the First Task on each Step.


\subsection{Step 1}

This automatic: here $\B_0=\mE_0=T(V)$ and $\gamma_0\colon\B_0\to \mE_0$ is $\gamma_0=\id$. So $\gamma_0(r)=r$, $r\in\Gc_0$, and thus
\[
\mE_1(\bs\lambda)=T(V)/\lg r-\lambda_r: r\in \Gc_0\rg.
\]

\subsection{Step $k+1$: How do we compute $\gamma_k(r)$?}\label{sec:task1}

Assume we have already computed $\mE_k(\bs\lambda)$. We look for $\gamma_k(r)\in\mE_k$, for each $r\in\Gc_k$.

Let us set, for short, $\gamma\ig\gamma_k$, $\mE\ig\mE_k$. We fix $r\in\Gc_k$.

\begin{taskone}
Compute $\gamma(r)\in \mE$ such that
\begin{align}\label{eqn:gamma1}
\rho(\gamma(r))=\gamma(r)\ot 1 + 1\ot r.
\end{align}
\end{taskone}

The solution, see \eqref{eqn:solution}, is given by the following recursive procedure:

\subsubsection{} We start by proposing
\[\gamma(r)\overset{?}{=} r.\]
If this satisfies \eqref{eqn:gamma1}, we are done. 

\subsubsection{} If not, we get an expression of the form:
\begin{align*}
\rho(r)=r\ot 1 + 1\ot r + \sum_{j=1}^{p(1)} t_j\ot t^j,
\end{align*}
where $t_j, t^j\in T(V)$ and $\deg(t_j)+\deg(t^j)<\gr(r)$, $1\leq j\leq p(1)\in\N$. 

We may assume, without lost of generality, that $\deg(t^j)< \deg(t^{j+1})$, $1\leq j\leq p(1)-1$. Notice that, as we only consider {\it strictly} increasing degrees, it may happen that $t^j$ is a {\it sum} of monomials $t^j=\sum_i m_{ij}$, each one of the with degree $\deg(t^j)$, that is $t^j$ does not necessarily represent a single monomial.

Moreover, as $t_{p(1)}\in\k$ by Lemma \ref{lem:tp} below, we may assume that $t_{p(1)}=1$.

Next, we propose:
\[\gamma(r)\overset{?}{=} r-t^{p(1)}.\]
If this satisfies \eqref{eqn:gamma1}, we are done. 

\subsubsection{} If not, notice that we have that
\begin{align*}
\rho(r-t^p)&=[r\ot 1 + 1\ot r + \sum_{j=1}^p t_j\ot t^j] - [t^p\ot 1 + 1\ot t^p + \sum_{j=1}^{p(2)} s_j\ot s^j]\\
&=(r-t^p)\ot 1 + 1\ot r + \sum_{j=1}^{p-1} t_j\ot t^j - t^p\ot 1 - \sum_{j=1}^{p(2)} s_j\ot s^j
\end{align*}
where $s_j, s^j\in T(V)$ and $\deg(s_j)+\deg(s^j)<\deg(t^p)$, $1\leq j\leq q$.

Hence, we can re-write this as
\begin{align*}
\rho(r-t^p)&=(r-t^p)\ot 1 + 1\ot r + \sum_{j=1}^{\ell} d_j\ot d^j
\end{align*}
where $d_j, d^j\in T(V)$ and $\deg(d_j)+\deg(d^j)<\deg(t^p)$, $1\leq j\leq q$, and $\deg d^j<\deg d^{j+1}$, $1\leq j\leq q-1$.

As $r-t^{p(1)}$ does not satisfy \eqref{eqn:gamma1}, we have that $\sum_{j=1}^{\ell} d_j\ot d^j\neq 0$. 

With the same argument as before, we may assume that $d_\ell=1$. We propose:
\[
\gamma_1(r)\overset{?}{=} r-t^p-d^\ell.
\]
and we repeat the previous analysis. 

\subsubsection{} It should be noticed that in each step we subtract to our candidate $\gamma(r)$ a term of lesser degree each time. This process finishes after a finite number of steps and we obtain a solution:
\begin{align}\label{eqn:solution}
\gamma(r)\overset{!}{=} r-t_1^{p_1}-t_2^{p_2}-\dots -t_m^{p_m}.
\end{align}

\subsubsection{} We complete this argument with the proof of our claim.
\begin{lem}\label{lem:tp}
$t_{p(1)}\in\k$.
\end{lem}
\pf
Set $p=p(1)$. Assume that $\deg(t_p)>0$. On the one hand, we have that
\begin{multline*}
(\id\ot\rho)\rho(r)=r\ot 1\ot 1+ 1\ot r\ot 1+1\ot 1\ot r \\
+\sum_j 1\ot t_j\ot t^j+\sum_j t_j\ot \rho(t^j).
\end{multline*}
On the other, this expression is equal to
\begin{multline*}
(\Delta\ot\id)\rho(r)=r\ot 1\ot 1+ 1\ot r\ot 1+1\ot 1\ot r +\sum_j \Delta(t_j)\ot t^j.
\end{multline*}
That is,
\begin{align}\label{eqn:igualdad-tp}
\sum_j 1\ot t_j\ot t^j+\sum_j t_j\ot \rho(t^j)=\sum_j \Delta(t_j)\ot t^j.
\end{align}
We write
\[
\rho(t^p)=t^p\ot 1 + 1\ot t^p + \sum_{i=1}^q s_i\ot s^i.
\]
for some $s_i,s^i\in T(V)$, with $\deg(s^i)<\deg(s^{i+1})$.
In particular, if we restrict on each side of the equality  \eqref{eqn:igualdad-tp} to those terms whose degree in the third tensorand is bigger (which is the case for $\deg(t^p)$), we have that, on the right hand side, these terms are
\[
\Delta(t_p)\ot t^p=t_p\ot 1\ot t^p + 1\ot t_p\ot t^p + \uDelta(t_p)\ot t^p.
\]
On the left hand side, these terms are
\[
1\ot t_p\ot t^p + t_p\ot \rho(t^p)= 1\ot t_p\ot t^p + \delta_{\deg(s^q),\deg(t^p)} t_p\ot s_q\ot s^q.
\]
That is,
\[
t_p\ot 1\ot t^p + \uDelta(t_p)\ot t^p=\delta_{\deg(s^q),\deg(t^p)} t_p\ot s_q\ot s^q.
\]
But this equality implies (by comparing the tensor of degree $\deg(t_p)$ in the first tensorand), that
$s^q=t^p$ and, more importantly, that
\[
\uDelta(t_p)=0.
\]
But this means that $t_p$ is a primitive element in $T(V)/\lg \Gc_0\rg$, with degree less than $\deg(r)$, which is a contradiction. Thus, $t_p\in \k$.
\epf

\subsection{The Hopf algebras}\label{sec:task2}

Assume we have already computed all 
$\gamma_k(r)$, $r\in\Gc_k$ so we have a new family of cleft objects $\mA_{k+1}=\mA_k/\lg\gamma_k(r)-\lambda_r\rg$. 

The cocycle deformation $\mL_{k+1}=L(\mA_{k+1},\mH_{k+1})$ is computed as a quotient of $\mL_{k}$, more precisely we have $\mL_{k+1}=\mL_k/\lg \tilde{r}-\lambda_r(1-g_r)\rg$, where 
\[
\tilde{r}=\delta(\gamma(r))-g_r\ot \gamma(r).
\]
That is we reach our Second Task.
\begin{tasktwo}
Compute $\delta(\gamma(r))$.
\end{tasktwo}
This is a routinary computation, see \ref{sec:remark-rho-delta}, (which can be, nevertheless, very much involved).

We shall provide  in \ref{sec:comp-alt} a series of tools to carry out these two tasks.

\subsection{The end} 
The reiteration of the two tasks described above ends with a family of Hopf algebra quotients of $T(V)\# H$, the algebras $\mL_{\ell+1}$, which is in turn a complete list of liftings of $V$.

It also provides a list of cleft objects for each lifting, as well as some families of pre-liftings, that is Hopf algebras that project onto the liftings in the list.

It is important to remark here that, in this diagonal case, it is always possible to define a stratification with {\it at most} four components $\Gc=\Gc_0\sqcup\Gc_1\sqcup\Gc_2\sqcup\Gc_3$; hence this strategy involves at most four steps (the first one being straightforward).

\subsection{Computational alternatives}\label{sec:comp-alt}

We see that the lifting problem relies in computing certain coaction formulae $\rho\colon\mA\to \mA\ot \mH$ or $\delta\colon\mA\to \mL\ot \mA$, both
induced by $\Delta:T(V)\to T(V)\ot T(V)$. 

\

We explain in detail two alternatives to perform these computations.

\subsubsection{A \texttt{GAP} algorithm}\label{sec:gap}

We explain how to solve the tasks described in \ref{sec:task1} and \ref{sec:task2}.

\begin{taskone}Find $\gamma(r)$.
\end{taskone}

Assume you want to find $\gamma(r)$ satisfying \eqref{eqn:gamma1} and you guess $\gamma(r)=r$. Then you need to compute $\rho(r)$ and check the equality.

In this case, $\rho:\mA_k\to \mA_k\ot \mH_k$; moreover we look for an expression of an element in $\mE_k$, that is we may restrict to 
$\rho:\mE_k\to \mE_k\#\k G\ot \B_k$, by \eqref{eqn:rho-i}, where $G=\lg g_i\rg_{i\in\I_n}\subseteq H$. Furthermore, for computational pourposes, we may consider $\hat{\rho}:\mE_k\to \mE_k\hat{\ot}\B_k$, where we denote $\hat{\ot}$ the tensor structure induced by the braiding in $\Rep G$. Finally, recall that $\rho(y_i)=y_i\ot 1+g_i\ot x_i$, $i\in\I_n$, so, in this setting 
\[
d_i\ig\hat{\rho}(y_i)=y_i\ot 1+1\ot x_i, \qquad i\in\I_n.
\]
In other words, to compute $\rho(r)$, where $r$ is a certain expression $r=\phi_r(y_1,\dots,y_n)$ involving sums of monomials on the letters $\{y_i\}_{i\in\I_n}$, we need to compute this same expression, now on the letters $\{d_i\}_{i\in\I_n}$, inside $\mE_k\hat{\ot}\B_k$. This is the algebra $F$ generated by $\{y_i,x_i\}_{i\in\I_n}$ with relations
\begin{align}
\label{rels-comm} x_iy_j=q_{ij} y_jx_i;&\\
\label{rels-E} \text{defining relations for }&\mE_k,  & \text{defining relations for }\mB_k.
\end{align}
We write this algebra in \texttt{GAP} (notice that we can obviate the symbol $\hat{\ot}$). 

To compute $\rho(r)$, we set 
\[d_i=y_i+x_i\in F, \ i\in\I_n\] which stands for the coaction $\rho(y_i)$ and compute $\phi_r(d_1,\dots,d_n)\in F$. 

Then, $\gamma(r)=r$ if and only if 
\begin{align}\label{eqn:check}
\phi_r(d_1,\dots,d_n)-\phi_r(y_1,\dots,y_n)-\phi_r(x_1,\dots,x_n)=0.
\end{align}

If we notice that $\gamma(r)\neq r$, then, as explained above, we take the highest order term $t^{p(1)}=t^{p(1)}(x_1,\dots, x_n)$ in \eqref{eqn:check}, set
$r'=r-t^{p(1)}$ and start over, which means we compute  $\phi_{r'}(d_1,\dots,d_n)=\phi_r(d_1,\dots,d_n)-t^{p(1)}(d_1,\dots, d_n)\in F$ and check if $\gamma(r')=r'$, that is if
\begin{align}\label{eqn:check-1}
\phi_{r'}(d_1,\dots,d_n)-\phi_{r'}(y_1,\dots,y_n)-\phi_r(x_1,\dots,x_n)=0.
\end{align}

If this is not the case, we start over... 

\begin{exa}
See \ref{sec:cleft5} for an explicit application of this idea.
\end{exa}

\begin{tasktwo}Compute $\delta(\gamma(r))$
\end{tasktwo}
The same ideas can be used to compute $\delta(\gamma(r))$, for the coaction $\delta:\mA_k\to \mL_k\ot \mA_k$ and thus give a presentation of the quotient Hopf algebra $\mL_{k+1}=L(\mA_{k+1},\mH_{k+1})$. 

Notice that $\gamma(r)\in \mE_k\subset \mA_k$, so it is an expression $\gamma(r)=\psi(y_1,\dots,y_n)$ in the variables $\{y_i\}_{i\in\I_n}$ and thus $\delta(\gamma(r))=\psi(\delta(y_1),\dots,\delta(y_n))$. Also, notice that $\gamma(y_i)=x_i$ and thus
\[
\delta(\gamma(x_i))=\delta(y_i)=a_i\ot 1+g_i\ot y_i, \qquad i\in\I_n.
\]
In this case we cannot obviate the group $G$, nor the generators $g_i$, $i\in\I_n$.

We write in \texttt{GAP} the algebra $E$ generated by $\{a_i,g_i,y_i\}_{i\in\I_n}$ with relations
\begin{align}
\label{rels-comm-L} g_ia_j&=q_{ij} a_jg_i, & y_ia_j&=a_jy_i, & g_iy_j&=y_jg_i;
\end{align}
\begin{align}
\label{rels-L} &\text{defining relations for }\mL_k, & \text{defining relations for }\mE_k.
\end{align}
Next, we set 
\[e_i\ig a_i+g_iy_i, \ i\in\I_n,\] 
which stands for the coaction $\delta(y_i)$ and compute $\psi(e_1,\dots,e_n)\in E$. 

In this case the answer we look for, namely the deformed relation $\tilde{r}$, is given in a single step by:
\[
\tilde{r}=\psi(e_1,\dots,e_n)-g_r\psi(y_1,\dots, y_n).
\] 
Here $g_r\in\lg g_1,\dots,g_n\rg\leq G$ is the group-like determined by $r$, see \eqref{eqn:gamma1}.

\begin{exa}
See \ref{sec:lift5} for an explicit application of this idea.
\end{exa}

\subsubsection{Explicit coproduct formulae}\label{sec:aar}

Recall that there is a distinguished set $\widetilde{\Gc}$ of relations in $\Gc$, given by the powers $x_\alpha^{N_\alpha}$, for $\alpha\in\mO(V)$ a root of Cartan type. 
Set $\Gc_0=\Gc\setminus \widetilde{\Gc}$. There is a more direct approach to computing the subfamily of liftings $L$ for which there is a projection $\widetilde{\B}(V)\ig T(V)/\lg \Gc_0\rg\twoheadrightarrow L$, that is the family of Hopf algebras $L$ that satisfy the relations $r=0$, $r\in\Gc_0$.
In particular, the set $\bs\Lambda$ can be described here as
\[
\{\bs\lambda=(\lambda_\alpha)_{\alpha\in\mO(V)}: \lambda_\alpha=0 \text{ if }  \chi_\alpha^{N_\alpha}\neq \eps\}
\]
and for each $\bs\lambda\in\Lambda$, we have 
 \begin{align*}
\mE(\bs\lambda)&\simeq \widetilde{\B}(V)/\lg \{ y_\alpha^{N_\alpha}-\lambda_\alpha |\alpha\in\mO(V)\},
 \end{align*}
 as algebras. 
 
The point is that, in this setting, we can derive explicit formulae for the coproducts $\Delta(x_\alpha^{N_\alpha})$, $\alpha\in\mO(V)$, for  $\Delta:\widetilde{\B}(V)\to \widetilde{\B}(V)\ot \widetilde{\B}(V)$. 

We need some notation. Let $\{\beta_1,\dots,\beta_M\}$ be an enumeration of the roots in $\mO$; so $M=|\mO|$. If ${\bf n}=(n_1,\dots,n_M)\in\Z_{\geq0}^M$, we set 
\begin{align*}
\underline{\bf n}&=n_1\beta_1 + \dots +n_M\beta_M, & x^{\bf n}&=x_{\beta_M}^{n_M}\dots x_{\beta_1}^{n_1}.
\end{align*}

Set $Z(V)=\k\{x_\alpha^{N_\alpha}:\alpha\in\mO\}$. Then this is a braided sub-Hopf algebra of $ \widetilde{\B}(V)$ by \cite[Theorem 4.13]{A-distinguished} and thus there are $r_{{\bf n},{\bf m}}(\alpha)\in\k$ such that
\begin{align}\label{eqn:coproduct}
\Delta(x_{\alpha}^{N_{\alpha}})=x_{\alpha}^{N_{\alpha}}\ot 1+ 1\ot x_{\alpha}^{N_{\alpha}}
+\sum_{\underline{\bf n}+\underline{\bf m}=N_\alpha \alpha }r_{{\bf n},{\bf m}}(\alpha) \ x^{\bf n}\ot x^{\bf m}.
\end{align}

Let $\pi_{\bs\lambda}\colon T(V)\twoheadrightarrow\mE(\bs\lambda)$ be the natural algebra projection, 
set $\lambda_\alpha({\bf m})=\pi_{\bs\lambda}(\bx^{\bf m})\in\k$ and consider the elements:
\begin{align*}
u_\alpha(\bs\mu)=\sum_{\underline{\bf n}+\underline{\bf m}=N_\alpha \alpha} \lambda_\alpha({\bf m})r_{{\bf n},{\bf m}}(\alpha)  
	a_{\beta_M}^{n_MN_{\beta_M}}\cdots a_{\beta_1}^{n_1N_{\beta_1}} g_{\beta_M}^{m_MN_{\beta_M}}\cdots g_{\beta_1}^{m_1N_{\beta_1}}.
\end{align*}

A large class of liftings can be classified in this way.
\begin{theorem}\cite[Theorem 3.1]{GJ}
For each $\bs\lambda\in\bs\Lambda$, the Hopf algebra $\ub(\bs\lambda)$ defined as the quotient of $T(V)\#H$ by the ideal generated by the relations
	\begin{align*}
r&=0, & r &\in\Gc_0;\\
\notag a_\alpha^{N_\alpha}&=\mu_\alpha(1-g_\alpha^{N_\alpha})-u_\alpha(\bs\mu), & \alpha&\in\mO(V),
	\end{align*}
is a lifting of $V$. Conversely, if $A$ is a Hopf algebra whose infinitesimal braiding is a principal realization $V\in\ydh$ and such that the relations $\Gc_0$ hold in $A$, then there is $\bs\lambda\in\bs\Lambda$ such that $A\simeq \ub(\bs\lambda)$.
\end{theorem}

In particular, this family of liftings is determined by the computation of the scalars $r_{{\bf n},{\bf m}}(\alpha)$. These can be computed in terms of the Lie algebra structure of a certain Lie algebra $\n$ associated to $(V,c)$: namely the exact sequence of braided Hopf algebras $Z(V)\stackrel{\iota}{\longrightarrow} \widetilde{\B}(V)\stackrel{\pi}{\longrightarrow}\B(V)$ gives raise to 
\begin{align}\label{eqn:sequence}
\B(V)\stackrel{\pi^*}{\longrightarrow}\L(V)\stackrel{\iota^*}{\longrightarrow}\z(V),
\end{align}
where $\L(V)$, resp. $\z(V)$, denote the graded dual of $\widetilde{\B}(V)$, resp. $Z(V)$. 

Let us assume for a moment that 
\begin{align}\label{eqn:hypothesis-AAR}
\hat{q}_{\alpha,\beta}^{N_\beta}=1, \quad \alpha,\beta\in\mO(V).
\end{align} 
Then $\z(V)\simeq U(\n)$, for $\n=\mP(\z(V))$. A basis of $\n$ is given by $\xi_\beta=\iota^*(y_\beta^{(N_\beta)})$, $\beta\in\mO(V)$, of the divided powers that generate $\L(V)$ and then
\begin{align}\label{eqn:c-t}
r_{{\bf n},{\bf m}}(\alpha)&=\prod_{i=1}^{k}\frac{1}{n_i!}\prod_{j=1}^{l}\frac{1}{m_j!}\,c(\alpha),
\end{align}
where $c(\alpha)\in\k$ is the unique scalar such that, in the algebra $U(\n)$,
\begin{align}\label{eqn:c-omega}
\xi_{\beta_M}\a{m_M}\dots \xi_{\beta_1}\a{m_1}\xi_{\beta_M}\a{n_M}\dots\xi_{\beta_1}^{n_1}=c(\alpha)\xi_\alpha + \text{other monomials}.
\end{align}
Recall that $\sum_{i=1}^M (n_i+m_i)\beta_i=N_\alpha\alpha$.

The general case, namely when condition \eqref{eqn:hypothesis-AAR} is removed, is obtained by observing that any braiding $(q_{ij})_{i,j\in\I}$ is twist equivalent to a braiding $(\hat{q}_{ij})_{i,j\in\I}$ satisfying this condition.

\begin{exa}
This idea is used in \S\ref{sec:example} to deal with the case $N\neq 5$.
\end{exa}

\section{Example}\label{sec:example}

We complete the classification of Hopf algebras of type $B_2$. When $N\neq 5$ is odd, the classification is given in \cite{BDR} and is recovered with our method above. When $N$ is even, the strategy is carried on by using ideas in \cite{AAR} to compute certain coproducts, as explained in \ref{sec:aar}. The case $N=5$, for a particular braiding matrix, cannot be solved with these ideas and the use of the computer seems to be imperative, as the deformation of the quantum Serre relations provokes a massive deformation of the other relations, leading to Hopf algebras with a very complicated presentation.

\subsection{Nichols algebras of Cartan type $B_2$}

Let $N\in\N_{\geq 3}$, $q\in\G_N'$. We set $M=\ord q^2$; i.e. $M=N$ if $N$ is odd, and $M=N/2$ if $N$ is even.

We consider here a matrix $\qb=(q_{ij})_{i,j\in\I_2}$ of Cartan type $B_2$; that is,
$q_{11}=q$, $q_{12}q_{21}=q^{-2}$, $q_{22}=q^2$.
The Nichols algebra $\B_{\qb}\coloneqq\B(V)$ is presented by generators $x_1$, $x_2$ and relations
\begin{align}\label{eq:qserre-B2}
&x_{1112}, & & & &x_{221},
\\ \label{eq:PRV-B2}
&x_{1}^N, & & x_{112}^M, & &x_{12}^N, & & x_2^M.
\end{align}
Let $\Bt_{\qb}$ be the algebra presented by generators $x_1$, $x_2$ and relations \eqref{eq:qserre-B2}, $\cHt=\Bt_{\qb}\# H$. In this Hopf algebra,
\begin{align}\label{eq:coprod-PRV-simple-roots}
\Delta(x_1^N)&=x_1^N\ot 1+ g_1^N\ot x_1^N; &
\Delta(x_2^M)&=x_2^M\ot 1+ g_2^M\ot x_2^M.
\end{align}
For the other powers of root vectors, the coproduct is computed in \cite{AAR}:

\begin{enumerate}[leftmargin=*,label=\rm{(\alph*)}]
	\item If $N$ is odd, then
	\begin{align}\label{eq:coproduct-PRV-odd-1}
	\Delta(x_{12}^N)=x_{12}^N\ot 1+ g_{12}^{N}\ot x_{12}^N
	+(1-q^{-2})^N q_{21}^{\frac{N(N-1)}{2}} x_1^N g_{2}^{N}\ot x_2^N; 
	\\ \label{eq:coproduct-PRV-odd-2}
	\begin{aligned}
	\Delta(x_{112}^N)=&x_{112}^N\ot 1+ g_{112}^{N}\ot x_{112}^N
	\\&
	+(1-q^{-1})^N(1-q^{-2})^N q_{21}^{N(N-1)} x_1^{2N}g_{2}^{N}\ot x_2^N \\&
	+2(1-q^{-1})^N(1+q)^N q_{11}^{\frac{N(N-1)}{2}}q_{21}^{\frac{N(N-1)}{2}} \, x_1^N g_{12}^{N}\ot x_{12}^N.
	\end{aligned}
	\end{align}
	\item If $N$ is even, then
	\begin{align}\label{eq:coproduct-PRV-even-1}
	\begin{aligned}
	\Delta(x_{12}^N)=&x_{12}^N\ot 1+ g_{12}^{N}\ot x_{12}^N
	+(1-q^{-2})^N q_{21}^{M(N-1)}\, x_1^N g_{2}^{N}\ot x_2^{N} \\&+ (1-q^{-2})^M q_{21}^{M^2} x_{112}^M g_{2}^{M}\ot x_2^M; 
	\end{aligned}
	\\ \label{eq:coproduct-PRV-even-2}
	\begin{aligned}
	\Delta(x_{112}^M)=&x_{112}^M\ot 1+ g_{112}^{N}\ot x_{112}^M 
	\\ &
	+(1-q^{-1})^M(1-q^{-2})^M q_{21}^{M(M-1)} x_1^{N} g_{2}^{M} \ot x_2^M.
	\end{aligned}
	\end{align}
\end{enumerate}

\begin{rem}\label{rem:q-serre-deformed-G5}
	If either $\chi_1^3\chi_2=\eps$ or else $\chi_1\chi_2^2=\eps$, then $q\in\G_5'$ and the matrix is 
	$\qb=\begin{pmatrix} q & q^2 \\ q & q^2 \end{pmatrix}$.
\end{rem}

\pf
Indeed, assume that $\chi_1^3\chi_2=\eps$. By evaluation in $g_1$,
$q_{11}^3q_{12}=1$, so $q_{12}=q^{-3}$ and then $q_{21}=q$. 
By evaluation in $g_1$, $1=q_{21}^3q_{22}=q^5$, so $q\in\G_5'$.

A similar computation solves the case $\chi_1\chi_2^2=\eps$.
\epf

\medbreak

To define the cleft objects and the corresponding liftings, we denote by
\begin{align}\label{eq:lambda-B2}
\bs\lambda &= (\lambda_{1}, \lambda_{2}, \lambda_{3}, \lambda_{4}, \lambda_{5}, \lambda_{6})
\end{align}
the corresponding family of scalars, appearing in the same order as \eqref{eq:qserre-B2}, \eqref{eq:PRV-B2}; that is, $\lambda_i:=\lambda_{r_i}$ is the scalar corresponding to the relation $r_i$.

\subsection{Case $N=2M$ even, $M>1$}
We fix the following stratification:
\begin{align}\label{eq:stratif-B2-even}
\Gc_0 &= \{x_{1112},x_{221}\}, & 
\Gc_1 &= \{x_{1}^{2M}, x_2^M \}, & 
\Gc_2 &= \{ x_{112}^M \}, & 
\Gc_3 &= \{x_{12}^{2M}\}.
\end{align}

\subsubsection{Cleft objects}
We look for the collection of algebras $\mA\in\Cleft(\mH)$; hence we have to compute $\mE=\mE(\bs\lambda)$, depending on $\bs\lambda=(\lambda_r)_{r\in \Gc}$, in such a way that $\mA=\mE\# H$.

By Remark \ref{rem:q-serre-deformed-G5}, $\chi_1^3\chi_2\neq\eps$,  $\chi_1\chi_2^2\neq \eps$. Hence $\lambda_{1}=\lambda_{2}=0$, and then
$\mE_1=\B_1=\Bt_{\qb}$.
Moreover, $\rho_1=\Delta:\mA_1\to \mA_1 \ot \mH_1$. By \eqref{eq:coprod-PRV-simple-roots},
\begin{align*}
\rho_1(y_1^{2m}) &= y_1^{2M} \ot 1 + g_1^{2M} \ot x_1^{2M},
&
\rho_1(y_2^{m}) &= y_2^{M} \ot 1 + g_2^{M} \ot x_2^{M},
\end{align*}
so $\mE_2=\mE_1/\langle y_1^{2M}-\lambda_{3}, y_2^M-\lambda_4 \rangle = T(V) /\langle y_{1112}, y_{221}, y_1^{2M}-\lambda_{3}, y_2^M-\lambda_4 \rangle$.

For the next step, we use \eqref{eq:coproduct-PRV-even-2} to prove that \begin{align*}
\rho_2(y_{112}^M)=&y_{112}^M\ot 1+ g_{112}^{N}\ot x_{112}^M,
\end{align*}
so $\mE_3 = T(V) /\langle y_{1112}, y_{221}, y_1^{2M}-\lambda_{3}, y_2^M-\lambda_4, y_{112}^M-\lambda_{5} \rangle$.

Finally, by \eqref{eq:coproduct-PRV-even-1} we have that
\begin{align*}
\rho_3(y_{12}^N)=&y_{12}^N\ot 1+ g_{12}^{N}\ot x_{12}^N.
\end{align*}
Hence we have that
\begin{align*}
\mE=\mE(\bs\lambda) = T(V) /\langle y_{1112}, y_{221}, y_1^{2M}-\lambda_{3}, y_2^M-\lambda_4, y_{112}^M-\lambda_{5}, y_{12}^N-\lambda_{6} \rangle.
\end{align*}

\subsubsection{Liftings}
Now we compute the liftings $\ub(\bs\lambda)$. According with the procedure above, we have to compute the Hopf algebras $\mL_k=\mL_k(\bs\lambda)$ such that $\mA_k$ is a $(\mL_k,\mH_k)$-biGalois object, $k=0,1,2,3$; so $\ub(\bs\lambda)=\mL_3(\bs\lambda)$.

We start with $\mL_0=T(V)\# H$. As $\mA_1=\mH_1$, we have that $\mL_1=\mH_1$ as well, and $\delta_1$ is just the comultiplication $\Delta:\mA_1\to \mL_1 \ot \mA_1$. Thus
\begin{align*}
\mL_2 &=\mL_1/\langle a_1^{2M}-\lambda_{3}(1-g_1^{2M}), a_2^M-\lambda_4(1-g_2^M) \rangle
\\
&=T(V)\# H/\langle a_{1112}, a_{221}, a_1^{2M}-\lambda_{3}(1-g_1^{2M}), a_2^M-\lambda_4(1-g_2^M) \rangle.
\end{align*}

Now we use \eqref{eq:coproduct-PRV-even-2} and the defining relations of $\mA_2$ to see that 
\begin{align*}
\delta_2(y_{112}^M)=& \big( a_{112}^M
+(1-q^{-1})^M(1-q^{-2})^M \lambda_4 \, a_1^{2M} g_{2}^{M}\big)
\ot 1+ g_{112}^{N}\ot y_{112}^M.
\end{align*}
Using the defining relations of $\mL_2$, we have that
\begin{align*}
\mL_3 &= \mL_2 /\langle a_{112}^M
+(1-q^{-1})^M(1-q^{-2})^M\lambda_4 \lambda_{3}(1-g_1^{2M}) g_{2}^{M}
-\lambda_{5} (1-g_{112}^M) \rangle.
\end{align*}

Finally, by \eqref{eq:coproduct-PRV-even-1} we have that
\begin{multline*}
\delta_3(y_{12}^N)=\Big( 
a_{12}^N +(1-q^{-2})^N \lambda_{4}^2\, a_1^N g_{2}^{N} + (1-q^{-2})^M \lambda_{4} a_{112}^M g_{2}^{M}  \Big) \ot 1\\
+ g_{12}^{N}\ot y_{12}^N.
\end{multline*}
Hence $\ub(\bs\lambda)$ is the quotient of $T(V)\# H$ by the relations
\begin{align*}
& a_{1112}=0, \qquad a_{221}=0, \qquad a_1^{2M}=\lambda_{3}(1-g_1^{2M}), \qquad a_2^M=\lambda_4(1-g_2^M),
\\
& a_{112}^M=\lambda_{5} (1-g_1^{2M}g_2^M)
-(1-q^{-1})^M(1-q^{-2})^M \lambda_{3}\lambda_4 \, (1-g_1^{2M}) g_{2}^{M},
\\
& a_{12}^N=\lambda_{6}(1-g_1^Ng_2^N)-(1-q^{-2})^{2M} \lambda_{3} \lambda_{4}^2 \, (1-g_1^{2M}) g_{2}^{2M} 
\\ & \qquad \qquad 
+ (1-q^{-2})^{2M} (1-q^{-1})^M \lambda_{3} \lambda_4^2 \, (1-g_1^{2M}) g_{2}^{2M}
\\ & \qquad \qquad 
-(1-q^{-2})^M \lambda_{4}\lambda_{5} (1-g_1^{2M}g_2^M)g_{2}^{M}.
\end{align*}

\subsection{Case $N=5$}\label{sec:b2-deformed}

In this part we assume $\qb=\begin{pmatrix}
	q&q^2\\q&q^2
	\end{pmatrix}$, $q\in\mathbb{G}_5$. This is the degenerate case. When $N=5$ but the braiding matrix is a different one, it can be treated with the ideas for $N$ odd in \ref{sec:N-odd}.

\medbreak

We fix the following stratification:
\begin{align}\label{eq:stratif-B2-order5}
\Gc_0 &= \{x_{1112},x_{221}\}, & 
\Gc_1 &= \{x_{1}^5, x_2^5 \}, & 
\Gc_2 &= \{x_{12}^5\}, & 
\Gc_3 &= \{x_{112}^5\}.
\end{align}
According to this choice, we shall consider parameters $\bs\lambda=(\lambda_1,\lambda_2,\lambda_3,\lambda_4,\lambda_5)$ associated to each generator. Observe that in this setting (and in this setting only), we see that we may have  
\[
\chi_1^3\chi_2=\chi_1\chi_2^2=\eps. 
\]
Hence $\lambda_{1}, \lambda_{2}$ may be nonzero scalars. In other words, the quantum Serre relations do not hold in any lifting; the cases in which these relation can be deformed give rise to very complicated expressions for the deformation of the powers of the root vectors. To deal with these cases, we use the computer program \texttt{GAP}, as explained in \S \ref{sec:gap}. We include in the Appendix a copy of the codes we use to compute the deformations.

We refer the reader to the \texttt{GAP} files as well as the \texttt{log} files storaged in the authors webpages:
\begin{enumerate} 
\item \texttt{B2-cleft.(g|log)} for the computation of the algebras $\mE(\bs\lambda)$.
\item \texttt{B2-lift.(g|log)} for the computation of the liftings $\ub(\bs\lambda)$.
\end{enumerate}
Both of these files have comments explaining the steps in the language of the present article.

\subsubsection{Cleft objects}\label{sec:cleft5}
We start with $\mE_0=T(V)$ and thus 
\[
\mE_1=T(V) /\langle y_{1112}-\lambda_1, y_{221}-\lambda_2\rangle.
\]
Now, we set $\mA_1=\mE_1\# H$, $\mH_1=\B_1\# H=T(V)\#H/\lg x_{1112}, x_{221}\rg$ 
and consider the right coaction $\rho_1:\mA_1\to \mA_1 \ot \mH_1$. 
By \eqref{eq:coprod-PRV-simple-roots}, we get
\begin{align*}
\rho_1(y_1^{5}) &= y_1^{5} \ot 1 + g_1^{5} \ot x_1^{5},
&
\rho_1(y_2^{5}) &= y_2^{5} \ot 1 + g_2^{5} \ot x_2^{5},
\end{align*}
so $\mE_2=\mE_1/\langle y_1^{5}-\lambda_{3}, y_2^5-\lambda_4 \rangle$.

Next we need to compute $\rho_2(y_{12})^5$, for $\rho_2:\mA_2\to \mA_2 \ot \mH_2$, where $\mA_2=\mE_2\#H$ and $\mH_2=\mH_1/\lg x_{1}^5, x_{2}^5\rg$. We use \texttt{GAP} and obtain
\[
\rho_2(y_{12})^5=y_{12}^5\ot 1+ g_{12}^{5}\ot x_{12}^5.
\]
Therefore, $\mE_3=\mE_3/\langle y_{12}^{5}-\lambda_{5}\rangle$. Finally, if $\mA_3=\mE_3\#H$ and $\mH_3=\mH_2/\lg x_{12}^5\rg$ and 
$\rho_3:\mA_3\to \mA_3 \ot \mH_3$ is the coaction, we use \texttt{GAP} again to see that
\[
\rho_3(y_{112})^5=y_{112}^5\ot 1+ g_{112}^{5}\ot x_{112}^5
\]
and hence  we have that
\begin{align*}
\mE(\bs\lambda) = T(V) /\langle y_{1112}-\lambda_1, y_{221}-\lambda_2, y_1^{5}-\lambda_{3}, y_2^5-\lambda_4, y_{112}^5-\lambda_{5}, y_{12}^5-\lambda_{6} \rangle.
\end{align*}

\subsubsection{Liftings}\label{sec:lift5}

We start with $\mL_0=T(V)\# H$ and
\begin{align*}
\mL_1 &=T(V)\# H/\langle a_{1112}-\lambda_1(1-g_{1112}^{5}), a_{221}-\lambda_2(1-g_{221}^{5}) \rangle.
\end{align*}
If $\delta_1\colon\mA_1\to \mL_1 \ot \mA_1$ is the left coaction, as the elements in $\Gc_1$ are primitive, we see that
\begin{align*}
\mL_2 &=\mL_1/\langle a_1^{5}-\lambda_{3}(1-g_1^{5}), a_2^5-\lambda_4(1-g_2^5) \rangle.
\end{align*}

Let us set $\delta_2\colon\mA_2\to \mL_2 \ot \mA_2$ is the left coaction. Observe that
\begin{align*}
\delta_2(y_1)&=a_1\ot 1+g_1\ot y_1, \qquad\qquad \delta_2(y_2)=a_2\ot 1+g_2\ot y_2,\\
\text{so }& \qquad  \delta_2(y_{12})=a_{12}\ot 1+g_1g_2\ot y_{12} + (1-q^{-2}) a_1g_2\ot y_2.
\end{align*}

We need to compute $\delta_2(y_{12})^5$. We use \texttt{GAP} and get
\begin{align*}
\delta_2(y_{12}^5)=& a_{12}^5\ot 1 + g_{12}^{5}\ot y_{12}^5 + s_{12}\ot 1
\end{align*}
where
\begin{align*}
s_{12}= & (-q-2q^2-3q^3-4q^4)\lambda_{1}\lambda_{2}  g_{2}^2g_{1}a_{2}a_{12}
+(2q+2q^2+q^4)\lambda_{1}\lambda_{2}^2          \\
&
+ (-10q-5q^2+5q^3+10q^4)\lambda_{3}\lambda_{4}   g_{2}^5g_{1}^5 \\
&
+ (10q+5q^2-5q^3-10q^4)\lambda_{3}\lambda_{4}    g_{2}^5\\
& 
+ (-2q^2-q^3-2q^4)\lambda_{1}\lambda_{2}^2       g_{2}^4g_{1}^2 
+ (-2q+q^3+q^4)\lambda_{1}\lambda_{2}^2           g_{2}^2g_{1}.
\end{align*}
Hence, 
\begin{align*}
\mL_3 &= \mL_2 /\langle a_{12}^5 + s_{12} - \lambda_5(1-g_{12}^5) \rangle.
\end{align*}

Finally, we consider the left coaction $\delta_3\colon\mA_3\to \mL_3 \ot \mA_3$, so that
\begin{multline*}
\delta_3(y_{112})= a_{112}\ot 1 + g_{112}\ot y_{112}\\
+ q(1-q^{-2}) a_{1}g_{12}\ot y_{12} 
+ (2+2q+q^{2})  a_{1}^2g_{2}\ot y_{2} 
\end{multline*}
and compute 
\begin{align*}
\delta_3(y_{112})^5=& a_{112}^5\ot 1 + g_{112}^{5}\ot y_{112}^5 +s_{112}\ot 1,
\end{align*}
where
\begin{align*}
s_{112}= &
(5q+5q^2+10q^3+5q^4)\lambda_{1}\lambda_{2}	g_{221}  		a_{112}^2a_{1}^2 \\
& +(5q+5q^2+5q^4)\lambda_{1}^2\lambda_{2}		g_{221}g_{1112}	a_{112}a_{1}\\
&
+(-5q-10q^2-10q^3-5q^4)\lambda_{1}^2\lambda_{2}	g_{221}		a_{112}a_{1}\\
&
+\big[(3q+q^2-q^3+2q^4)\lambda_{1}^3\lambda_{2}  +(15q+5q^2+20q^3+10q^4)\lambda_{1}\lambda_{2}^2\lambda_{3}  \\
&\qquad+(25q-25q^2-25q^3+25q^4)\lambda_{3}^2\lambda_{4}
\\
&\qquad\qquad +(-10q+20q^2-20q^3+10q^4)\lambda_{3}\lambda_{5}	\big]			g_{2}^5g_{1}^{10} \\
&
+\big[ (-15q-5q^2-20q^3-10q^4)\lambda_{1}\lambda_{2}^2\lambda_{3}\\
&\qquad+(-50q+50q^2+50q^3-50q^4)\lambda_{3}^2\lambda_{4}\\
&\qquad\qquad  +(10q-20q^2+20q^3-10q^4)\lambda_{3}\lambda_{5}	\big]	
g_{2}^5g_{1}^5 \\
&
+ (-8q-6q^2-4q^3-2q^4)\lambda_{1}^3\lambda_{2}				g_{2}^3g_{1}^4 \\
&+ (2q-q^2+q^3-2q^4)\lambda_{1}^3\lambda_{2}			g_{2}^4g_{1}^7 \\
&
+ (25q-25q^2-25q^3+25q^4)\lambda_{3}^2\lambda_{4}				g_{2}^5 \\
&
+ (3q+6q^2+4q^3+2q^4)\lambda_{1}^3\lambda_{2}				g_{2}^2g_{1}. 
\end{align*}
Hence $\ub(\bs\lambda)$ is the quotient of $T(V)\# H$ by the relations
\begin{align*}
a_{1112}&=\lambda_1(1-g_{1112}), \qquad a_{221}=\lambda_2(1-g_{221}), \\
a_1^{5}&=\lambda_{3}(1-g_1^{5}), \qquad\qquad a_2^5=\lambda_4(1-g_2^5),\\
a_{12}^5&=\lambda_{5}(1-g_{12}^5)-s_{12},\\
a_{112}^5&=\lambda_{6} (1-g_{112}^5)-s_{112}.
\end{align*}

\subsection{Case $N\neq 5$ odd}\label{sec:N-odd}

In this part, we recover the classification in \cite{BDR}. We also allow the case $N=5$ but $\qb\neq \begin{pmatrix}
	q&q^2\\q&q^2
	\end{pmatrix}$, $q\in\mathbb{G}_5'$, which is the degenerate case solved in \ref{sec:b2-deformed}.

We consider the filtration as in \eqref{eq:stratif-B2-order5}. By Remark \ref{rem:q-serre-deformed-G5}, $\chi_1^3\chi_2\neq\eps$,  $\chi_1\chi_2^2\neq \eps$. Hence $\lambda_{1}=\lambda_{2}=0$, so $\mA_1=\mH_1=\mL_1=\B_{\qb}\# H$.

Now we use the formulas \eqref{eq:coprod-PRV-simple-roots}, \eqref{eq:coproduct-PRV-odd-1} and \eqref{eq:coproduct-PRV-odd-2} to prove that 
\begin{align*}
\mE(\bs\lambda) &= T(V) /\langle y_{1112}, y_{221}, y_1^N-\lambda_{3}, y_2^N-\lambda_4, y_{12}^N-\lambda_{5}, y_{112}^N-\lambda_{6} \rangle,
\end{align*}
and $\ub(\bs\lambda)$ is the quotient 
of $T(V)\# H$ by the relations
\begin{align*}
& a_{1112}=0, \qquad a_{221}=0, \qquad a_1^N=\lambda_{3}(1-g_1^N), \qquad a_2^N=\lambda_4(1-g_2^N),
\\
& a_{12}^N= \lambda_{5} (1-g_1^Ng_2^N)-(1-q^{-2})^N \lambda_{3}\lambda_4 (1-g_1^N) g_{2}^{N}
,
\\
& a_{112}^N=\lambda_{6}(1-g_1^{2N}g_2^N) -(1-q^{-1})^N(1-q^{-2})^N \lambda_{3}^2 \lambda_4  (1-g_1^N)^2g_{2}^{N} 
\\ & \qquad\qquad  -2(1-q^{-1})^N(1+q)^N \lambda_{3} \lambda_5 (1-g_1^N) g_{12}^{N}.
\end{align*}

\section{A more general framework}\label{sec:gral}

We move now to the more general the setup described in \ref{sec:more-general}:
\begin{enumerate}
	\item A cosemisimple Hopf algebra $H$.
	\item A braided vector space $(V,c)$ such that the ideal $\mJ(V)$ defining the Nichols algebra $\B(V)$ is finitely generated.
\end{enumerate}
That is, we allow infinite-dimensional Nichols algebras $\B(V)$, non-diagonal braidings and non-principal realizations $V\in \ydh$.

\subsection{Infinite-dimensional Nichols algebras}
Jordan and super Jordan braidings naturally appear when we study Yetter-Drinfeld modules over abelian groups which are not of diagonal type. Indeed these are the unique examples of indecomposable Yetter-Drinfeld modules whose Nichols algebra has finite Gelfand Kirillov dimension \cite{AAH-triang}. 

The braided vector spaces are denoted $\mathcal V(2,\eps)$, with $\epsilon=1$ for Jordan braiding and $\eps=-1$ for super Jordan braiding. They are of dimension two, with a basis $\{x_1,x_2\}$ such that
\begin{align}\label{eq:braiding-jordan-superJ}
c(x_i \ot x_1) &= \eps x_1 \ot x_i, &
c(x_i \ot x_2) &= (\eps x_2 + x_1) \ot x_i, &
i &=1,2.
\end{align}
Both Nichols algebras $\B(\mathcal V(2,\eps))$ have GK dimension two, and the defining relations are:
\begin{align*}
 x_2x_1-x_1x_2-\frac{1}{2}x_1^2&=0, & &\text{when }\eps=1;
\\
 x_1^2=0, \quad
 x_2x_{21}-x_{21}x_2-x_1x_{21}&=0, & &\text{when }\eps=-1,
\end{align*}
where $x_{21}\coloneqq(\ad_{c} x_2)x_1 =x_2x_1+x_1x_2$.

Let $G$ be an abelian group. A realization $\mathcal V(2,\eps)\in\ydg$, $\eps=1,-1$, is given by a  Jordanian, respective super Jordanian, YD-triple
$\mathcal D=(g,\chi,\eta)$, for some $g\in G$, $\chi\in\widehat{G}$ and a $(\chi,\chi)$-derivation $\eta\colon \k G\to \k$.

Our strategy to compute liftings also applies for these examples. For Jordan braidings, the procedure is easy since we have a unique relation. 
For super Jordan braidings, we have a stratification of two steps since $x_1^2$ is primitive in $T(\mathcal V(2,-1))$ while the other relation is not primitive. Applying the strategy, we obtain \cite[Propositions 4.2 \& 4.4]{AAH-liftings}: Every pointed Hopf algebra with coradical $\Bbbk G$ and Jordan, respectively super Jordan, infinitesimal braiding is of the form
$\mathfrak U(\mathcal D,\lambda)$, where $\mathcal D=(g,\chi,\eta)$ is a Jordanian, respective super Jordanian, YD-triple and $\lambda\in\Bbbk$ is subject to the condition
\begin{align*}
\lambda &=0 & \mbox{if }\chi^2 \neq \eps.
\end{align*}
Here, $\mathfrak U(\mathcal D,\lambda)$ is the quotient of $T(\mathcal V(2,\eps))\# \Bbbk G$ by the relations
\begin{align*}
 x_2x_1-x_1x_2-\frac{1}{2}x_1^2 &= \lambda(1-g^2), & &\text{when }\eps=1,
\\ x_1^2= \lambda(1-g^2), \quad  x_2x_{21}-x_{21}x_2-x_1x_{21}&=-\lambda (2x_2 + x_1 g^2), & &\text{when }\eps=-1.
\end{align*}

\subsection{Non abelian groups} When the group $G$ is not abelian, the classification of all $V\in\ydg$ with $\dim\B(V)<\infty$ is not complete. However, when such $V$ is known and a presentation of $\B(V)$ is given, the strategy also applies and has been used to construct liftings in \cite{GV1} and \cite{GV2}, which are, a fortiori, cocycle deformations of $\B(V)\# \k G$. In this case, the braiding is not diagonal: rather it is determined by a conjugacy class $\mO\subset G$ in such a way that there is a basis $\{x_i:i\in\mO\}$ of $V$ for which
\begin{equation}\label{eqn:rack-br}
c(x_i\ot x_j)=q_{ij}\,x_{iji^{-1}}\ot x_i, \qquad i,j\in\mO,
\end{equation}
and $(q_{ij})_{i,j\in\mO}$ is a family of scalars satisfying certain conditions. More generally, if $X$ is a rack and $q:X^2\to\k$ is a rack 2-cocycle, then these data determine a braided vector space $V(X,q)$. If $G$ is such that $V\in\ydg$, the strategy can be applied, see loc.cit.

\subsubsection{A warning}\label{sec:warning}
In this setting, we may not be a able to make the assumption on \eqref{eqn:assume}; namely that each $\Gc_k=\{r\}$ for a single $r$. Rather, we shall have a sub-object $\Gc_k\subset \B_k$ in $\ydh$. The procedure remains the same, but an extra task is needed: in the language of \ref{sec:quotients}, if $S_k\coloneqq\k\{rg_r^{-1}:r\in\Gc_k\}$, we need to check that the comodule map $\varphi_{\bs\lambda}\colon S_k\to \mA_k$ extends to an algebra morphism $Y_k\coloneqq \k\lg S_k\rg\to \mA_k$. 

\begin{exa}\label{exa:s3-pointed}
Let $\mO_{(12)}\subset G=\s_3$ be the conjugacy class of the transposition $(12)$ and consider the braiding on $V=\k\{x_{(12)}, x_{(13)}, x_{(23)}\}$ as in \eqref{eqn:rack-br} with $q_{ij}=-1$, all $i,j\in\mO$. This defines the unique $V\in\ydg$ with $\dim\B(V)<\infty$; the realization is given by 
\[
x_i\mapsto i\ot x_i, \quad g\cdot x_i=\sgn(g)\,x_{gig^{-1}}, \qquad i\in\mO_{(12)}, g\in \s_3;
\]
and $\B(V)$ is the algebra with generators $x_{(12)},x_{(13)},x_{(23)}$ and relations
\begin{align}\label{eqn:rels-s3}
\begin{split}
x_{(12)}^2=x_{(13)}^2=x_{(23)}^2&=0,\\
x_{(12)}x_{(13)}+x_{(23)}x_{(12)}+x_{(13)}x_{(23)}&=0,\\
x_{(13)}x_{(12)}+x_{(12)}x_{(23)}+x_{(23)}x_{(13)}&=0. 
\end{split}
\end{align}
The liftings $\ub(\lambda)$, $\lambda\in\k$, are defined as the quotients of $T(V)\#\k\s_3$ modulo:
\begin{align*}
x_{(12)}^2=x_{(13)}^2=x_{(23)}^2&=0,\\
x_{(12)}x_{(13)}+x_{(23)}x_{(12)}+x_{(13)}x_{(23)}&=\lambda(1-(132)), \\
x_{(13)}x_{(12)}+x_{(12)}x_{(23)}+x_{(23)}x_{(13)}&=\lambda(1-(123)). 
\end{align*}
These Hopf algebras arise via the strategy by considering the cleft objects $\mA(\lambda)=\mE(\lambda)\#\k\s_3$, where the algebra $\mE(\lambda)$ is generated by $x_{(12)},x_{(13)},x_{(23)}$, with relations
\begin{align*}
x_{(12)}^2=x_{(13)}^2=x_{(23)}^2&=0,\\
x_{(12)}x_{(13)}+x_{(23)}x_{(12)}+x_{(13)}x_{(23)}&=\lambda,\\
x_{(13)}x_{(12)}+x_{(12)}x_{(23)}+x_{(23)}x_{(13)}&=\lambda. 
\end{align*}
In  particular, all liftings are cocycle deformations of $\B(V)\#\k \s_3$.
\end{exa}

\subsection{Copointed Hopf algebras}\label{sec:copointed}

A Hopf algebra $H$ is called copointed if $H_0=\k^G$, for some non-abelian group $G$. When $G$ is finite, there is a braided equivalence of categories $F\colon \ydg\to\ydgdual$, thus any $V\in\ydg$ as above gives rise to $F(V)\in\ydgdual$ with $\dim\B(F(V))<\infty$. 

In particular, if $(V,c)$ is a braided vector space with a principal realization $V\in\ydg$, then it determines a realization $V\in\ydgdual$ and the same ideas can be used to compute the liftings of $V\in\ydgdual$. 
This is the path used in \cite{GV2} to complete the classification of copointed Hopf algebras over $\s_4$, and to recover the classification in \cite{AV}, over $\s_3$. As well, the liftings over $\k^{\Dm_m}$ in \cite{FGM} can also be recovered with this method.

\begin{exa}\label{exa:s3-copointed}
Consider once again $\mO_{(12)}\subset G=\s_3$ and $V$ as in Example \ref{exa:s3-pointed}. The realization $V\in\ydgdual$ is given by
\[
x_i\mapsto \sum_{g\in G}\sgn(g) \delta_g\ot x_{g^{-1}ig}, \quad \delta_g\cdot x_i=\delta_{g,i}\,x_i, \qquad i\in\mO_{(12)}, g\in \s_3.
\]
Here $\{\delta_g\}_{g\in \s_3}$ denotes the standard basis of idempotents in $\k^{\s_3}$.
In particular, this is not a principal realization. However, the strategy applies. The Nichols algebra $\B(V)$ has the same presentation as in \eqref{eqn:rels-s3} and the liftings are presented as quotients of $T(V)\#\k^{\s_3}$ with the relations
\begin{align*}
x_{(13)}^2&=  (\lambda_1 - \lambda_2 )(\delta_{(12)} + \delta_{(123)} ) + \lambda_1 (\delta_{(23)} + \delta_{(132)} ),\\
x_{(23)}^2&=  \lambda_2 (\delta_{(13)} + \delta_{(123)} ) + (\lambda_2 - \lambda_1 )(\delta_{(12)} + \delta_{(132)} ),\\
x_{(12)}^2&= - \lambda_1 (\delta_{(23)} + \delta_{(123)} ) - \lambda_2 (\delta_{(13)} + \delta_{(132)} ),\\
x_{(12)}x_{(13)}&+x_{(23)}x_{(12)}+x_{(13)}x_{(23)}=0, \\
x_{(13)}x_{(12)}&+x_{(12)}x_{(23)}+x_{(23)}x_{(13)}=0. 
\end{align*}
These Hopf algebras were defined in \cite[Definition 3.4]{AV}  and can also be obtained using the strategy; the associated cleft objects are of the form $\mA(\lambda_1,\lambda_2)=\mE(\lambda_1,\lambda_2)\#\k^{\s_3}$, where $\mE$ is generated by $x_{(12)},x_{(13)},x_{(23)}$ with relations
\begin{align*}
x_{(12)}^2=\lambda_1, \qquad x_{(13)}^2=\lambda_2, \qquad x_{(23)}^2&=-\lambda_1-\lambda_2,\\
x_{(12)}x_{(13)}+x_{(23)}x_{(12)}+x_{(13)}x_{(23)}&=0,\\
x_{(13)}x_{(12)}+x_{(12)}x_{(23)}+x_{(23)}x_{(13)}&=0. 
\end{align*}
In  particular, all liftings are cocycle deformations of $\B(V)\#\k^{\s_3}$.
\end{exa}

\subsection{Non-principal realizations}\label{sec:non-principal}

As we saw in \ref{sec:copointed}, the strategy also applies for certain non-principal realizations $V\in\ydgdual$. However, these realizations are in the image of a principal realization under a category equivalence $\ydg\to \ydgdual$. 

With full generality, it can also be applied for any non-principal realization $V\in\ydh$. We recall from \cite[4.1]{GV2} that such realization of a braided vector space $(V,c)$ with basis $\{x_i\}_{i\in\I}$ and braiding $c(x_i\ot x_j)=\sum_{k,l\in\I}c_{kl}^{ij}\, x_k\ot x_l$ is determined by 
matrix coefficients $\{\mu_{ij}\}_{i,j\in\I}\subset H^\ast$ and comatrix elements $\{e_{ij}\}_{i,j\in\I}\subset  H$ describing the $H$-action and coaction:
\begin{align*}
h\cdot x_i&=\sum_{j\in\I}\mu_{ij}(h) x_j, & x_i&\mapsto \sum_{j\in I}e_{ij}\ot x_j, \qquad i\in \I, \, h\in H;
\end{align*}
subject to 
\begin{align*}
\sum_{j\in\I}\mu_{ij}(h_{(1)})e_{jk}h_{(2)}&=
\sum_{j\in\I}\mu_{jk}(h_{(2)})h_{(1)}e_{ij}, &\text{and} & & \mu_{jl}(e_{ik})=c_{lk}^{ij},
\end{align*}
for all $h\in H$ and each $i,j,k,l\in\I$, which express the compatibility in $\ydh$ and the coherence of the braidings, respectively.

%
%
%
%
%
%


\begin{thebibliography}{AAG1}
\bibitem[A+]{AAGMV} {\sc Andruskiewitsch, N.}, {\sc Angiono, I.}, {\sc Garc\'ia Iglesias, A.}, {\sc Masuoka, M.},
{\sc Vay, C.} \emph{Lifting via cocycle deformation}. J. Pure Appl. Alg. {\bf 218} (4), 684--703 (2014).
%
\bibitem[AAG]{AAG}  {\sc Andruskiewitsch, N.}, {\sc Angiono, I.}, {\sc Garc\'ia Iglesias, A.}
\emph{ Liftings of Nichols algebras of diagonal type I. Cartan type A},  Int. Math. Res. Not. IMRN {\bf 9}, 2793--2884 (2017).

\bibitem[AAH1]{AAH-triang} 
{\sc Andruskiewitsch, N.}, {\sc Angiono, I.}, {\sc Heckenberger, I.}
\emph{On finite GK-dimensional Nichols algebras over abelian groups}. \texttt{arXiv:1606.02521}. Mem. Amer. Math. Soc., to appear.


\bibitem[AAH2]{AAH-liftings} \bysame 
\emph{Liftings of Jordan and super Jordan planes}. \texttt{arXiv:1512.09271}. Proc. Edinb. Math. Soc., II. Ser., to appear.


\bibitem[AAR]{AAR}  {\sc Andruskiewitsch, N., Angiono, I., Rossi Bertone, F.} \emph{The quantum divided power algebra of a finite-dimensional Nichols algebra of diagonal type}. Math. Research Letters \textbf{24} 619--643 (2017).

\bibitem[AC]{AC}  {\sc Andruskiewitsch, N., Cuadra, J.} \emph{On the structure of (co-Frobenius) Hopf algebras}. J.
Noncommut. Geom. \textbf{7}  83--104 (2013).
%

\bibitem[AG]{AG} {\sc Andruskiewitsch, N.}, {\sc Garc\'ia Iglesias, A.}, 
\emph{Twisting Hopf algebras from cocycle deformations}, Annali dell'Università di Ferrara {\bf 63} 221--247, (2017).

%
%
%
%
%
%
%
%
\bibitem[AS1]{AS-lift-meth} {\sc Andruskiewitsch, N.}, {\sc Schneider, H.J.}, \emph{Lifting of quantum linear spaces and pointed Hopf algebras of order $p^3$}. J. Algebra {\bf 209}, 658--691 (1998).

\bibitem[AS2]{AS-pointed}  \bysame,
\emph{Pointed Hopf algebras}, ``New directions in Hopf algebras'',
MSRI series Cambridge Univ. Press; 1--68 (2002).
%
\bibitem[AS3]{AS-annals} \bysame,
\emph{On the classification of finite-dimensional pointed Hopf
algebras}, Ann. of Math. \textbf{171}, 375--417 (2010).
%
%
\bibitem[AV]{AV} {\sc Andruskiewitsch, N.}, {\sc Vay, C.}, \emph{Finite dimensional
Hopf algebras over the dual group algebra of the symmetric group in three letters},
Comm. Algebra {\bf 39}, 4507--4517 (2011).
%

\bibitem[A1]{A-jems} {\sc Angiono, I.} {\em A presentation by generators and relations of Nichols algebras
    of diagonal type and convex orders on root systems},  J. Europ. Math. Soc. {\bf 17} (2015) 2643--2671.

\bibitem[A2]{Ang-crelle} \bysame
{\em On Nichols algebras of diagonal type}, J. Reine Angew. Math.  \textbf{683} 189--251 (2013).

\bibitem[A3]{A-distinguished} \bysame \emph{Distinguished Pre-Nichols algebras}, Transf. Groups {\bf 21}, 1--33, (2016).

\bibitem[AnG]{AnG} {\sc Angiono, I.,  Garc\'ia Iglesias, A.}
\emph{Liftings of Nichols algebras of diagonal type II. All liftings are cocycle deformations}. \texttt{arXiv:1605.03113}.


\bibitem[BDR]{BDR} {\sc Beattie, M., D\u{a}sc\u{a}lescu, S., Raianu, S.} \emph{Lifting of
Nichols algebras of type $B_2$}. Israel J. Math. {\bf 132}  1--28 (2002).

%
%
%
\bibitem[FGM]{FGM} \textsc{Fantino, F.}, {\sc Garc\'ia, G.A.}, {\sc Mastnak, M.}
\textsl{On finite-dimensional copointed Hopf algebras over dihedral groups}, \texttt{arXiv:1608.06167}.

\bibitem[GJ]{GJ} {\sc Garc\'ia Iglesias, A.}, {\sc Jury Giraldi, J. M.}, 
\textsl{Liftings of Nichols algebras of diagonal type III. Cartan type $G_2$}, J. Algebra {\bf 478}, 506--568 (2017).

%
%
%
%
%
\bibitem[GV1]{GV1} {\sc Garc\'ia Iglesias, A.}, {\sc Vay, C.}
\textsl{Finite-dimensional pointed or copointed Hopf algebras over affine racks}, J. Algebra {\bf 397}  379--406 (2014).

\bibitem[GV2]{GV2} \bysame,
\textsl{Copointed Hopf algebras over $\s_4$}, J. Pure Appl. Alg. {\bf 222} (9), 2784--2809 (2018).

%
%
%

\bibitem[GM]{GM} {\sc Grunenfelder L., Mastnak, M.} \textsl{Pointed and copointed Hopf algebras as cocycle deformations}, \texttt{arxiv:0709.0120v2}.

\bibitem[Gu]{Gunther} {\sc Gunther, R.} \textsl{Crossed products for pointed Hopf algebras},
Comm. Algebra, {\bf 27}, 4389--4410 (1999).
\bibitem[H]{H} {\sc Heckenberger, I.} \emph{Classification of arithmetic root systems},
Adv. Math. \textbf{220},  59--124 (2009).
\bibitem[He]{Helbig} {\sc Helbig, M.} \emph{On the lifting of Nichols algebras}, Comm.
in Alg. \textbf{40} (2012) 3317--3351.
%
\bibitem[M]{M} {\sc Masuoka, A.} \textsl{Abelian and non-abelian second cohomologies of
quantized enveloping algebras}. J. Algebra, \textbf{320}, 1--47 (2008).
%
%
%

\bibitem[S]{S} {\sc Schauenburg, P.}, {\it Hopf bi-Galois extensions}, Comm. Algebra
\textbf{24}, 3797--3825, (1996).

%
\end{thebibliography}
\end{document}